\documentclass[11pt]{amsart}
\usepackage{amsmath,amsthm, amscd, amssymb, amsfonts}
\usepackage{t1enc}
\usepackage{latexsym}
\usepackage{stmaryrd}

\usepackage{verbatim}
\usepackage{amsgen,amsmath,amstext,amsbsy,amsopn,amsfonts,amssymb}

\newtheorem{theorem}{Theorem}[section]
\newtheorem{lemma}[theorem]{Lemma}
\newtheorem{corollary}[theorem]{Corollary}

\newtheorem{proposition}[theorem]{Proposition}

\theoremstyle{definition}

\theoremstyle{remark}
\newtheorem{remark}[theorem]{Remark}

\newcommand{\n}{{\mathfrak n}}
\newcommand{\vo}{{\mathfrak v}}
\newcommand{\wo}{{\mathfrak w}}

\newcommand{\ub}{{\mathfrak u}}

\newcommand{\g}{{\mathfrak g}}

\newcommand{\h}{{\mathfrak h}}

\newcommand{\z}{{\mathfrak z}}

\newcommand\Aut{\operatorname{Aut}}
\renewcommand{\ker}{\mbox{\rm Ker\,}}

\newcommand{\R}{{\mathbb R}}
\newcommand{\Hc}{{\mathbb H}}
\newcommand{\Oc}{{\mathbb O}}

\newcommand{\C}{{\mathbb C}}
\newcommand{\Lag}{{\mathcal L}}

\renewcommand\mod{\operatorname{mod}}

\numberwithin{equation}{section}

\newcommand{\F}{\mathbb F}

\begin{document}

\title{Horizontal submanifolds of groups of Heisenberg type}

\author[Kaplan, Levstein, Saal and Tiraboschi]
{A. Kaplan, F. Levstein, L. Saal and A. Tiraboschi}
\address{\noindent
CIEM (CONICET) and FaMAF (UNC), Ciudad Universitaria, C\'ordoba
5000, Argentina} \email{\{kaplan, levstein, saal,
tirabo\}@mate.uncor.edu}
\address{\noindent Kaplan also at: University of
Massachusetts, Amherst, MA 01003, USA}
\email{kaplan@mate.uncor.edu}
\thanks{\noindent This work was supported by CONICET, Antorchas,  FONCyT
 and Secyt (UNC)}

\begin{abstract}
We study maximal horizontal subgroups of Carnot groups of
Heisenberg type. We classify those of dimension half of that of
the canonical distribution (``lagrangians'') and illustrate some
notable ones of small dimension. An infinitesimal classification
of the arbitrary maximal horizontal submanifolds follows as a
consequence.
\end{abstract}

\maketitle

\begin{section}{Introduction }

 A general Carnot manifold is, by definition, endowed with a
 bracket-generating
 distribution. A horizontal submanifold is one whose tangent spaces lie in
 the distribution.

 Because the distribution is ``outvolutive'',
 horizontal submanifolds of large dimension are rare. At the same time, when
 one exists,
 there may be a continuum of others through each
 point, even sharing tangent
 and higher jet spaces there. Maximal ones are most natural to
 study on general grounds  and  appear in Geometric Control Theory
as  jet spaces of maps, or limits of minimal submanifolds, for
example.
  The goal of this article is to describe a representative class of maximal
horizontal submanifolds in a representative class of Carnot
manifolds.

Carnot groups are representative of Carnot manifolds in a strict
sense, the latter carrying canonical sheaves of the former. In
such a group, its Lie subgroups form a representative class of
submanifolds, and we will prove that any horizontal submanifold is
osculated everywhere by translates of horizontal subgroups. In
particular, the  possible dimensions are the same.

Let $N$ be a Carnot group ($N$ is for nilpotent). In terms of the
intrinsic grading of the Lie algebra
 $$\n={\rm Lie}(N)= \n_1\oplus ... \oplus \n_s,$$
 a horizontal subgroup is of the form $\exp(\ub)$, where $\ub\subset
 \n_1$ is a subspace satisfying $[\ub,\ub]=0$. One is so lead to describe
 the maximal abelian
 subalgebras contained in the generating subspace $\n_1$. But not much else can be said
 in  general, because there are just too many Carnot groups,
 even 2-step ones.

A natural subclass to study is that of groups of Heisenberg type.
Their role in subriemannian geometry has been compared to that of
euclidean spaces, or of symmetric spaces, in riemannian geometry
and the  analysis of elliptic operators \cite{CGN1}. In practice,
they remain the only Carnot groups where basic convexity questions
can be answered at all, or abundant domains fit for the Dirichlet
problem for its sublaplacian  can be actually constructed
\cite{CGN1} \cite{CGN2} \cite{Ka1}. In any case, they are an
obvious starting point. The fact that they have been a source of
unexpected examples in ordinary riemannian geometry and analysis
gives additional support for the choice
\cite{Ka2}\cite{L}\cite{DR}\cite{GV}\cite{DGGW}\cite{GW}\cite{Sz1}\cite{Sz2}

Let then $N$ be of Heisenberg type and use the standard notation
for the  grading of $\n$:
$$\n_1 = \vo, \qquad \n_2=\z=\ {\rm center}(\n)$$
as well as for the dimensions
$$n= \dim \vo,\qquad m=\dim \z.$$
For emphasis, the intrinsic distribution on $N$ has dimension $n$
and codimension $m$ (there will be no specific notation for the
dimension of $N$ itself, to avoid confusion).

We find that the dimension of a maximal horizontal subspace of
$\vo$ and, therefore, that of any maximal horizontal submanifold
of $N$, must be among the numbers
$$ {n\over 2},\ {n\over 3},\ ...\ ,\ {n\over  m+1 },$$
so the set ${\rm Hor}(N)$ of maximal horizontal subgroups can have
at most $m$ strata. In the ordinary Heisenberg groups, $m=1$,
$n=2k$, and all the maximal horizontal subgroups are
$k$-dimensional (``planes''). These are the Lagrangian subspaces
(maximal isotropic) of the obvious symplectic forms defined by the
bracket, which justifies much of the terminology used here.

In this paper we first describe ${\rm Hor}(N)$  in some notable
examples, enough to illustrate generic features, like non-trivial
stratifications, as well as peculiar ones, like 8-dimensional
distributions with no horizontal submanifolds of dimension $>1$.
In the second part we fully describe the first stratum,
$${\rm Lag}(N) = \{ U\in {\rm Hor}(N):\  \dim U = {n\over 2}\}$$
This is a real-analytic variety, which is sometimes empty,
sometimes it is a Lie group, and always is a finite union of
orbits of the analogous ``symplectic'' group $\Aut_o(\n)$,
consisting of the automorphisms of $\n$ that fix the  central
elements.

To be more specific and as we recall in the first section, we can
pick
$$\z=\R^m$$
 for any $m\geq 1$, and let $\vo$ be
any finite-dimensional module over the Clifford algebra $C(m)$.
This will be of the form
$$\vo=(\vo_m )^{p }$$
if $m\neq 3$ modulo 4, or
$$\vo=(\vo_m^+)^{p_+} \oplus
(\vo_m^-)^{p_-}$$ otherwise, where $\vo_m,\vo_m^{\pm},$ are
irreducible (the real spinor spaces). Hence, $\n$ is determined by
the integers $m,n,$ or $m, p_+, p_-$. The following is a table
showing the values that yield a non-empty ${\rm Lag}(N)$, together
with their structure under $\Aut_o(\n)$.

\begin{table}[h]

\begin{tabular}{|cccc|}
\hline
$m \,(\mod 8)$       & $p$ or $(p_+,p_-)$& $\sharp$ orbits & Lag($G$)  \\
      0               & any $p$           &   $p+1$ & $ \bigcup_r
 O(p)/O(r)\times O(p-r)$ \\
      1               &   any $p$         &    1           & $U(p)/O(p)$  \\
      2               &  any $p$           & 1         &  $U(p,\Hc)/U(p)$ \\
      3               &  any $(p,p)$       & 1         & $U(p,\Hc)$ \\
      4               &  any $p$           & $p+1$ & $\bigcup_r U(p,\Hc) /
 U(r,\Hc)\times U(p-r,\Hc)$ \\
      5               &  $p$ even        & 1      & $U(p)/U(p/2,\Hc)$  \\
      6               &  $p$ even        & 1      & $O(p)/O(p/2) \times
 O(p/2)$  \\
      7               &  any $(p,p)$     & 1      & $O(p)$  \\
       \hline
\end{tabular}

\vskip .5cm
 \caption{The variety of Lagrangians subspaces} \label{tablaLang}

\end{table}

In the last section we give explicit descriptions of ${\rm
Lag}(N)$ in terms of Plücker coordinates.

{\bf Acknowledgement:} We thank L. Capogna and N. Garofalo for
their advise and W. Dal Lago for his hospitality and patience.

\end{section}

\begin{section}{Lie algebras of Heisenberg type}

Let $V$ be an $\R $-vector space and let $B: V \times V \to \R $
be a non-degenerate symmetric bilinear form. Recall that a  {\it
Clifford algebra} associated to $(V,B)$ is a a pair $(C(V,B),
\theta)$, where $ C(V,B)$ is an $\R $-algebra, $\theta: V \to
C(V,B)$ is a linear function such that $\theta(x)^2= B(x,x)1$ for
each $x \in V$ and $(C(V,B), \theta)$ satisfies the following
universal property: if $(A, \mu)$ is a pair such that $A$ is an
$\R $-algebra and $\mu:V \to A$ is  linear and satisfies
$\mu(x)^2=B(x,x)1$, then there exists an algebra morphism
$\mu':C(V,B)\to A$ such that $\mu'\theta=\mu$ and $\mu'$ is unique
with respect to this property.

The Clifford algebra $(C(V,B),\theta)$ exists for every $(V,B)$
and can be obtained as a quotient of the tensor algebra $T(V)$ by
the ideal generated by $x\otimes x -B(x,x)1$. Moreover,
$(C(V,B),\theta)$ is unique modulo isomorphism of $\R $-algebras
and $V$ is naturally embedded in $C(V,B)$.

Let $C(m)$ denote $C(\Bbb R^m,-\langle x,y \rangle)$, where
$\langle \,,\,\rangle$ is the standard inner product on $\R^m$.
Every module over the algebra $C(m)$ is unitary and is the direct
sum of irreducibles modules. Up to isomorphism, there is precisely
one irreducible module $\vo_m$ over $C(m)$ for $m \not\equiv 3
\,(\mod 4)$, and two, $\vo_m^+,\vo_m^-$, for $m\equiv 3 \,(\mod
4)$ \cite{Hu}\cite{BtD}.

If $\vo$ is a $C(m)$-module with compatible inner product $(u,v)$,
we put a Lie algebra structure on
$$\n = \vo \oplus \R^m$$
by declaring $\R^m$ to be the center and  defining
$$[\ ,\ ]:\vo\wedge\vo \rightarrow \R^m$$
by
$$\langle z,[u,v] \rangle = (J_zu,v)$$
where $J_z$ denotes the Clifford action. This satisfies
$J_z^2=-\langle z,z \rangle I_{z}$, hence $[u,v] $ is indeed
skew-symmetric; Jacobi's identity is trivially satisfied.

A Lie algebra of Heisenberg type with center $\z=\R^m$ is then of
the form
$$\n =(\vo_m)^p \oplus \z,$$
for $m\not \equiv 3 \,(\mod 4)$ or
$$\n = (\vo_m^+)^{p_+} \oplus (\vo_m^-)^{p_-} \oplus \z$$
for $m \equiv 3 \,(\mod 4)$. These are mutually non-isomorphic
except that $(p_+,p_-)$ and $(p_-,p_+)$ give isomorphic Lie
algebras. One says that $\n$ is {\em irreducible} if the
corresponding Clifford module is irreducible.

The structure of $C(m)$ as associative algebra and the dimension
of its irreducible representations $\vo$ are listed in the table
below. Let $\R$, $\C$ and $\Hc$ denote the real, complex and
quaternionic numbers, and  $\R_k$, $\C_k$, $\Hc_k$ the algebra of
matrices of order $k$ with coefficients in $\R$, $\C$, $\Hc$,
respectively.

\begin{table}[h]

\begin{tabular}{|ccccc|}
\hline $m \,(\mod 8)$ & 0     &1      &2        &3  \\
\hline $C(m)$  \qquad &  \quad $\R_{2^k}$  \quad &\quad$\C_{2^k}$
\quad
 &\quad$\Hc_{2^{k-1}}$   \quad     &\quad$\Hc_{2^{k-1}}\otimes
 \Hc_{2^{k-1}}$  \quad \\
$\dim_\R\vo$    &\,\,\,  $2^{k}$     &\,\,\, $2^{k+1}$ &\,\,\,
 $2^{k+1}$        &\,\,\, $2^{k+1}$ \\
\hline
\end{tabular}

\vskip .5cm

\begin{tabular}{|ccccc|}
\hline $m \,(\mod 8)$   &4 &5 &6 &7\\
\hline $C(m)$  \qquad &\quad$\Hc_{2^{k-1}}$ \quad      &\quad
$\C_{2^k}$
 \quad     &\quad $\R_{2^k}$ \quad &\quad $\R_{2^k} \oplus \R_{2^k} $
 \quad
 \\
$\dim_\R\vo$    &\,\,\, $2^{k+1}$&\,\,\, $2^{k+1}$&\,\,\,
$2^{k}$&\,\,\,
 $2^{k}$\\
\hline
\end{tabular}

\vskip .5cm
 \caption{$C(m)$ and the dimension of irreducible representations. Here,
 $m=2k$ if $m$ is even,
$m=2k+1$ otherwise} \label{tabla1}

\end{table}

Let $C^+(m)$ denote the even Clifford algebra, generated by the
elements of even order $zz'$, and $z_1,\ldots,z_m$ an orthonormal
basis of $\z = \R^m$. Then the map
\begin{equation}
z_j \mapsto z_jz_m, \qquad \text{ for } 1 \le j \le m, \qquad
\text{ and } z_m \mapsto z_m
\end{equation}
extends to an automorphism of $C(m)$, which takes $C(m-1)$ onto
$C^+(m)$.

The element
$$K_m = J_{z_1} \ldots J_{z_m}$$
commutes with $C^+(m)$ and, when $m$ is odd,  with all of $C(m)$,
and

-- If $m \equiv 1,2 \, (\mod 4)$, $K_m^2=-1$ and $K_m^t = -K_m$.

-- If $m \equiv 0,3 \, (\mod 4)$, $K_m^2=1$ and  $K_m^t = K_m$.

-- If  $m \equiv 3 \,(\mod 4)$, $K_m$ acts on an irreducible
module as $\pm Id$. \noindent In the last case we denote by
$\vo_m^\pm$  the eigenspace of $K_m$ of eigenvalue $\pm 1$.

Let $\Aut(\n)$ be the group of Lie algebra automorphisms of $\n$,
$\Aut_o(\n)$ the subgroup  of elements acting trivially in the
center. Let $End_{C^+(m)}(\vo)$ denote the algebra of linear maps
on $\vo$ which commute with the action of $C^+(m)$. Then \cite{S}
\begin{equation} \label{eq:2.3}
{\rm Aut}_o(\n) = \{\xi \in End_{C^+(m)}(\vo): \xi^t J_z \xi =
J_z, \text{ for some } z \in \z, z \not=0 \}.
\end{equation}

\vskip .3cm

The algebras with $m=1,2,3,4,7,8$ can be described in terms of the
classical real division algebras $\R, \C, \Hc , \Oc,$ as follows.
Let $\F$ denote one of these  and consider  $\F^k$ as a real
vector space.

 To obtain those with $m=1,2,4,8$, take $\z=\F$ and $\vo= \F^p \times
 \F^p$.
 Then the bracket
$[\,,\,]: \vo \times \vo \to \z$ is
\begin{equation}\label{prodint4}
[({\bf x},{\bf y}),({\bf x}',{\bf y}')] = \sum_{j=1}^p x_j{y'_j} -
{x_j'}y_j
\end{equation}
and the  Clifford action is
\begin{equation}\label{action4}
J_z({\bf x},{\bf y}) = (-z \overline{{\bf y}},\overline{{\bf
x}}z), \qquad (z\in \F, \,\, {\bf x},{\bf y} \in \F^p).
\end{equation}

To obtain those with $m=1,3,7$, take $\F=\C, \Hc , \Oc,$
respectively, let $\z = \Im(\F)$ and $\vo = \F^p \times \F^q$.
Then the bracket $[\,,\,]: \vo \times \vo \to \z$ is
\begin{equation}\label{prodint3}
[({\bf x},{\bf y}),({\bf x}',{\bf y}')] = -\Im(\sum_{j=1}^p
x_j\overline{x'_j} + \sum_{k=1}^q \overline{y_k'}y_k ).
\end{equation}
and the Clifford action is
\begin{equation} \label{action3}
J_z({\bf x},{\bf y}) = (z{\bf x},{\bf y}z), \qquad (z\in \Im(\F),
\,\, {\bf x} \in \F^p,\,\, {\bf y} \in \F^q).
\end{equation}

Finally, note that
\begin{equation}\label{corchete3}
\langle ({\bf x},{\bf y}),({\bf x}',{\bf y}')\rangle =
\Re(\sum_{j=1}^k x_j\overline{x'_j}).
\end{equation}
is the natural inner product in $\F^k$.

\end{section}

\begin{section}{Horizontal submanifolds and subgroups}

\begin{proposition}\label{A} Let $S\hookrightarrow G$ be a horizontal
 submanifold of a
Carnot group $G$ and $g\in S$. Then there exist a unique Lie
subgroup $H\subset G$  such that $T_g(S)= T_g(g H)$. The subgroup
$H$ is horizontal and abelian.
\end{proposition}
\begin{proof} Since the distribution $\vo$ is left-invariant,
it is enough to assume that $S$ passes through the identity $e$
and that $g=e$. Let $\h \subset \g$ be the subspace spanned by the
left-invariant vector fields whose value at $e$ is tangent to $S$.
Let $X,Y\in \h$ and extend $X_e,Y_e$ to vector fields $\tilde{X}$,
$\tilde{Y}$, in a neighborhood of $e$, so that they are tangent to
$S$ along $S$. Therefore $[\tilde{X}, \tilde{Y}]$ will also have
this property; in fact, it may be assumed that $[\tilde{X},
\tilde{Y}]=0$ along $S$. Assuming also that $X,Y$ are linearly
independent, complete them to a basis $X_1=X$, $X_2=Y$, $X_3$, ...
,  of $\g_1$. Let $\{T_{\alpha j}\}$ be a basis of $\g_\alpha$ for
$\alpha\geq 2$ and write
$$\tilde{X} =\sum_i f_i  X_i  + \sum_{\alpha,j} \phi_{\alpha j}  T^{\alpha j}
 ,\qquad
 \tilde{Y} =\sum_i
g_i  X_i + \sum_{\alpha,j} \psi_{\alpha j} T_{\alpha j}$$ with
smooth coefficients. Now compute
$$0=[\tilde{X},\tilde{Y}] =  \sum_{i<j} (f_ig_j- f_jg_i)[X_i,X_j] +
\sum_{\alpha,j}((\tilde{X}\psi_{\alpha j})-(\tilde{Y}\phi_{\alpha
j}))T_{\alpha j} + \sum_{\alpha,\beta,j,k} \phi_{\alpha j}
\psi_{\beta k} [T_{\alpha j}, T_{\beta k}]$$ By horizontality, the
functions $\phi_{\alpha j}$ and $\psi_{\beta k}$  are constantly
equal to zero on $S$. Since $\tilde{X}$ and $\tilde{Y}$ are
tangent to $S$, the second and third sums vanish on $S$ and,
therefore,
$$ 0=  \sum_{i<j} (f_i(s)g_j(s)- f_j(s)g_i(s))
[X_i,X_j](s)$$  for $s\in S$. Evaluating at $s=e$ and recalling
that
 $f_1(e)=1$, $f_i(e) =0\ \forall i\not=1$, $g_2(e)=1$,
 $g_i(e) =0\ \forall i\not=2$,  we obtain
$[X,Y](e)=0$. By left-invariance, $[X,Y]=0$. We conclude that $\h$
is an abelian subalgebra, contained in $\vo$. Taking $H=
\exp(\h)$, the assertion follows.

\end{proof}

\begin{proposition} The maximal abelian subgroups of a 2-step Carnot group
$G$ are those of the form $U\cdot Z(G)$, where $U$ is a maximal
horizontal subgroup and $Z(G)$ is the center of $G$.
\end{proposition}
\begin{proof}
A maximal abelian subalgebra, as well as the center, are
automatically graded.
\end{proof}

\begin{proposition}\label{B} Let
${n}$ the dimension of the canonical distribution $\vo$ and $m$
its codimension. Then the dimension
 of any maximal horizontal submanifold of $G$ must be among the numbers $
 {n}/2,\ {n}/3,\ ... , {n}/(m+1)$.
\end{proposition}

\begin{proof} By Proposition \ref{A}, it is enough to prove the assertion for a
maximal horizontal subgroup, of the form $H=\exp (\h)$ with
$[\h,\h]=0$. In terms of the operators $J_z$, the commutativity is
expressed by $(J_z\h, \h)=0$, i.e., $J_\z\h\subset \h^\perp$.
Conversely, let $v\in \h^\perp$. Then $0=(v,J_Z\h)=(z, [v,\h])$
$\forall z$ and, therefore, $[v,\h]=0$. Since $\h$ is maximal
horizontal, $v\in \h$. Consequently, a subspace $\h\subset \vo$ is
maximal abelian if and only if
 $$J_\z (\h) = \h^\perp \label{X}$$
or, equivalently,
$$ \vo = \h \oplus (\sum J_{z_i}(\h)).$$
It follows that $2h\leq {n}\leq h + m h= h(m+1)$, so that
$$ {n}/(m+1)\leq h \leq {n}/2
$$
\end{proof}

Since  ${n}\geq 2^{{m-1\over 2}},$ Proposition \ref{B} is quite
restrictive and shows that for large $m$, the possible dimensions
are all close to largest possible one. The two extreme cases
correspond, respectively, to
$${n}/2:\qquad\qquad \ \ \h^\perp= J_{z}(\h)$$
for all $z\not=0$, and
$${n}/(m+1):\qquad\qquad  \h^\perp = \oplus_{i=1}^m J_{z_i}(\h)$$
for a basis of $\z$.

We next illustrate some possibilities that can occur. First,
consider  the lowest-dimensional  Heisenberg groups associated to
the division algebras as described at the end of last section.
Then
$$
{\rm Lie}(N_\C)=\C\oplus \Im(\C),\qquad {\rm Lie}(N_\Hc)=\Hc\oplus
\Im(\Hc),\qquad {\rm Lie}(N_\Oc)=\Oc\oplus \Im(\Oc).
$$
These are the only groups of Heisenberg type satisfying $n=m+1$.
Indeed, if $\{z_i\}$ is a basis of $\z$ and $v\in \vo$ is
non-zero, $v, J_{z_1}v,...,J_{z_m}v,$ is a basis of $\vo$ such
that $[v, J_{z_i}v]=z_i$. Therefore $\R v$ is maximal abelian in
$\vo$. A similar argument shows that they are the only irreducible
groups of Heisenberg type for which the dimension ${n}/(m+1)$ is
actually realized.

We conclude that the horizontal submanifolds of these groups are
all one-dimensional. The varieties of maximal horizontal subgroups
are
$${\rm Hor}(N_\C)= {\rm Lag}(N_\C)\simeq S^1$$
$${\rm Hor}(N_\Hc) \simeq \R P^3,\qquad {\rm Lag}(N_\Hc)=\emptyset$$
$${\rm Hor}(N_\Oc) \simeq \R P^7,\qquad {\rm Lag}(N_\Oc)=\emptyset$$

Finally, consider the case $m=8$ and $\n$ irreducible:
$$\n= \vo_8\oplus \R^8= (\Oc\times
\Oc)\oplus  \Oc.$$ Here $n=16$ and, therefore, the dimensions
allowed by Proposition \ref{B} are $2,4$ and $8$. In the next
section we will see that there are only two lagrangians
$$
{\rm Lag}(\n)= \{\wo_+ ,\wo_-\}.$$ Here we will see  that there is
none of dimension four and a $32$-parameter family of
two-dimensional ones; at the end,
$$ {\rm Hor}(\n)\cong \big(Gr_\R(2,8)\times
 \R P^7\times \R_+\big)\cup\{\wo_+ ,\wo_-\}.$$
When the parameter $t\in \R_+$ goes to $0$ or $\infty$, the limits
of the corresponding 2-dimensional subspaces lie, respectively, in
$\wo_+$ and $\wo_-$, so the terms in parenthesis define a natural
stratification of ${\rm Hor}(\n)$.

 Decompose $\vo_8 =
\wo_+ \oplus \wo_-$ as $C^+(8)$-module. Let $w \in \vo$, $w = v +
u$ with $v \in \wo_+, u \in \wo_-$. It easy to see that there
exists $z \in \z$ such that $w = v + J_zv$.

\begin{lemma} \label{propcaso8}  Let  $0 \not= v \in \wo_+$ and $0 \not= z=
 \in \z$.
\begin{enumerate}
\item The centralizer of $v +J_zv$ is  $J_{\z}(v+J_z^{-1}v)$.
\item Let $u,u' \in \z$ such that $u,u',z$ are linear independent,
then $J_u(v+J_z^{-1}v)$ and $J_{u'}(v+J_z^{-1}v)$ do not commute.
\end{enumerate}
\end{lemma}
\begin{proof}
Define $z'=z/\langle z,z \rangle$ and $t=\langle z,z \rangle$,
thus $J_z= tJ_{z'}$ and $J_z^{-1}= -t^{-1}J_{z'}$.

(1) It is clear that $v +J_zv$ and  $J_z(v+J_z^{-1}v) = J_zv + v $
commute. Now let $u \in z^\perp$ and $u'' \in \z$.  Then
\begin{align*}
\langle u'', [v +J_zv, J_u(v+J_z^{-1}v)] \rangle &= \langle u'',
[v , J_uv] \rangle + \langle u'', [tJ_{z'}v,
J_u(-t^{-1})J_{z'}^{-1}v]=
 \rangle \\
&= \langle u'',
[v , J_uv] \rangle - \langle u'', [J_{z'}v, J_uJ_{z'}^{-1}v] \rangle \\
&= \langle
J_{u''}v , J_uv \rangle - \langle J_{u''}J_{z'}v, J_uJ_{z'}v \rangle \\
& = \langle J_{u''}v , J_uv \rangle -\langle J_{u''}v , J_uv
\rangle =0.
\end{align*}
Thus, the centralizer of $v +J_zv$ contains $J_\z(v+J_z^{-1}v)$.
Recall that $ad(u) : \vo_8 \to \z$ is surjective for all $0\not=u
\in \vo_8$. Therefore, $\dim(\ker ad(v +J_zv))=8$, thus the
centralizer of $v +J_zv$ is equal to $J_\z(v+J_z^{-1}v)$.

(2) We can consider $u,u',z$ mutually orthogonally and $\langle
u,u \rangle= \langle u',u' \rangle =1$. Let $u'' \in {\rm
span}_\R\{u,u',z\}^\perp$, then
\begin{align*}
\langle u'', [J_u(v +J_z^{-1}v), J_{u'}(v+J_z^{-1}v)] \rangle &=
\langle u'', [J_u v , J_{u'}J_z^{-1} v] \rangle + \langle u'',
[J_uJ_z^{-1}v,=
 J_{u'}v] \rangle \\
&= -t\langle J_{u''}J_u v , J_{u'}J_{z'} v \rangle - t\langle
J_{u''}J_uJ_{z'}v, J_{u'}v=
 \rangle \\
&= -t\langle J_{u''} v , J_uJ_{u'}J_{z'} v \rangle - t\langle
J_{u''}v, J_{z'}J_uJ_{u'}v=
 \rangle \\
&= -2t \langle J_{u''}v, J_{z'}J_uJ_{u'}v \rangle.
\end{align*}
Thus, if we suppose that $J_u(v+J_z^{-1}v)$ and
$J_{u'}(v+J_z^{-1}v)$ commute we have $\langle J_{u''}v,
J_{z'}J_uJ_{u'}v \rangle = 0$. Since, we also have that $\langle
J_{w}v, J_{z'}J_uJ_{u'}v \rangle = 0$ for all $w \in {\rm
span}_\R\{u,u',z\}$. Then, $J_{z'}J_uJ_{u'}v=0$, a contradiction.
This proves the Lemma.
\end{proof}

The lemma implies that the 2-planes spanned by pairs $ u+J_zv,
v+J_zu,$ are all abelian and maximal with this property. To prove
that these are all those of dimension $<8$, let $u \in \vo$ and
write $u=(x,y)$ with respect to the decomposition
 $\vo = \wo_+ \oplus \wo_-$. Let $\Lag$ and $\Lag_1$
maximal isotropic subspaces of $\vo$ of dimension 2. $\Lag$
($\Lag_1$ respectively) has a basis $(x,y),(x',y')$ (resp.
$(x_1,y_1),(x_1',y_1')$) such that $x\perp x'$ (resp. $x_1 \perp
x_1'$) and $0=[(x,y),(x',y')]= xy'-x'y$ (resp. $0=
x_1y_1'-x_1'y_1$). If $\Lag=\Lag_1$, there must be $a,b,c,d \in
\R$ such that
\begin{alignat}2
x_1 &= ax +bx',& \qquad y_1 &= ay +by' \label{r1} \\
x_1' &= cx +dx',& y_1' &= cy + dy'. \label{r2}
\end{alignat}
Because of $x_1\perp x'_1$, it follows $\begin{pmatrix}a & c \\ b
& d  \end{pmatrix} \in O(2)$. Conversely, if $\begin{pmatrix}a & c
\\ b & d  \end{pmatrix} \in O(2)$ and $x_1,x'_1,y_1,y'_1$ are
defined as in formulas (\ref{r1}) and (\ref{r2}), then
$(x_1,y_1),(x_1',y_1')$ is a basis of $\Lag$ such that $x_1,x_1'$
is an orthonormal set. Thus, every maximal isotropic subspace
other than the $\wo_\pm$, is determined by a 2 dimensional
subspace of $\R^8$ (generated by $x,x'$) and a non-zero $y \in
\R^8$.

\end{section}

\begin{section}{Generalities on Lagrangians}

Let $N$ be a group of Heisenberg type with Lie algebra
 $\n=\vo\oplus\z $.
The maximal horizontal submanifolds of dimension $\dim \vo /2$ are
called Lagrangians. Via the exponential map, they correspond to
$[\ ,\ ]$-isotropic subspaces of $\vo$ of half the dimension,
which will also be called Lagrangians. They form a closed
real-analytic variety, which we  denote by ${\rm Lag}(N)$, ${\rm
Lag}(\n)$, or even ${\rm Lag}(\vo)$. In this section we describe
three general properties needed later: the relation with
$C^+(m)$-submodules, the appearance of the periodicity modulo 8
and the natural action of ${\rm Aut}(\n)$ on ${\rm Lag}(\n)$.

\begin{proposition}\label{prop:lperp} Let $\Lag\in {\rm Lag}(\vo)$.
Then  $\Lag^\perp\in {\rm Lag}(\vo)$,  $\vo = \Lag \oplus
\Lag^\perp$  and $J_z(\Lag) = \Lag^\perp$ for all non-zero $z\in
\z$.
\end{proposition}
\begin{proof} For $z \in \z$ and $x,y \in
\Lag$, $\langle J_z x, y \rangle = \langle z , [x,y] \rangle = 0$.
Thus $J_z(\Lag) = \Lag^\perp$ for any $z \in \z$. The other
results follow easily.
\end{proof}

It follows that double products $J_zJ_w$  preserves $\Lag$, hence

\begin{corollary}\label{cmasmodulo} Any Lagrangian $\Lag$  is a
$C^+(m)$-module.
\end{corollary}

Sometimes, to stress the dimension $m$ of the center $\z$,  it
will be useful to denote it by $\z_m$, so that $\z_m=\R^m$ and
$\z_m\times\z_r = \R^{m+r}$.

\begin{proposition} \label{prop:tensorial} Let $\vo_r$ be a
$C(r)$-module and let $\vo_m$  be a $C(m)$-module with $m  \equiv
0\ (\mod 4)$.  Then $\vo_m \otimes \vo_r$ is a $C(m+r)$-module,
with Clifford action
$$J_{(z,w)}:= J_z \otimes Id + K_m
\otimes J_w,$$ $(z,w) \in \z_m \times \z_r.$ The corresponding
algebra of Heisenberg type is
$$(\vo_m \otimes \vo_r)
\oplus (\z_m \times \z_r)$$ with  bracket
\begin{equation}\label{corchetetensorial}
[x\otimes u, y \otimes v] = (\langle u, v \rangle[x,y], \langle
K_mx,y \rangle[u,v]),
\end{equation}
$x,y \in \vo_m$, $u,v \in \vo_r$.
\end{proposition}
\begin{proof}
(\ref{corchetetensorial}) is checked   by taking inner product
with $z \in \z_m \times \z_r$ on both sides of the equation. The
rest of the assertions follow easily.
\end{proof}

\begin{corollary}\label{tensorirreducible} If $\vo_r$ is an irreducible
 $C(r)$-module, then
$\vo_8 \otimes \vo_r$ is an irreducible $C(8+r)$-module. Indeed,
$(\vo_8)^{\otimes s} \otimes \vo_r$ is an irreducible
$C(8s+r)$-module and
$$K_{8s+r} = (K_8^{r+1})^{\otimes s} \otimes
K_r.$$
\end{corollary}
\begin{proof}
The dimension of $\vo_8 \otimes \vo_r$ is $16 \times \dim(\vo_r)$
and we conclude the result from Table \ref{tabla1}.
\end{proof}

We can now see how the periodicity modulo 8 typical of Clifford
modules is reflected on these Lagrangians.

\begin{theorem}\label{cor3.5} If $\Lag_r\in {\rm Lag}(\vo_r)$ and
$\Lag_8\in {\rm Lag}(\vo_8)$, then
$$\Lag_8
\otimes \Lag_r + {\Lag_8}^\perp \otimes \Lag_r^\perp\ \in\ {\rm
Lag}(\vo_8 \otimes \vo_r).$$
\end{theorem}
\begin{proof} $\Lag^\perp$ is
Lagrangian from Proposition \ref{prop:lperp}. If $(x,u),(y,v) \in
\Lag' \otimes \Lag$, then $[x\otimes u, y \otimes v] =
[x,y]\langle u, w \rangle + \langle K_8x,y \rangle[u,v] = 0$,
because of $ \Lag$ and $\Lag'$ are Lagrangians. As $\Lag'^\perp$
and  $\Lag^\perp$ are isotropic, we have $[x\otimes u, y \otimes
v] = 0$ for $(x,u),(y,v) \in \Lag'^\perp \otimes \Lag^\perp$.

Now, let $(x,u) \in \Lag' \otimes \Lag$ and $(y,v) \in \Lag'^\perp
\otimes \Lag^\perp$. By Corollary \ref{cmasmodulo} we have that
$\Lag'$ is a $C^+(8)$-module, so $K_8x \in \Lag'$, then $[x\otimes
u, y \otimes v] = [x,y]\langle u, w \rangle + \langle K_8x,y
\rangle[u,v] =0$.
\end{proof}

Recall ${\rm Aut}_o(\n)$, the group of orthogonal automorphisms of
$\n$ that act trivially on $\z$. Then,

\begin{theorem} \label{transitively}
If $m \not\equiv 0  \,(\mod 4)$, $\Aut_o(\n)$ acts transitively on
${\rm Lag}(\n)$. If $m \equiv 0  \,(\mod 4)$ and $\n = (\vo_m)^p
\oplus \z$, then ${\rm Lag}(\n)$ is the union of $p+1$ orbits of
$\Aut_o(\n)$.
\end{theorem}
\begin{proof}
Let $\Lag, \Lag'\in {\rm Lag}(\n)$. By Corollary \ref{cmasmodulo},
$\Lag$ and $\Lag'$ are $C^+(m)$-modules of the same dimension.

If $m \not\equiv 0  \,(\mod 4)$, there must be an isomorphism
$\psi: \Lag \to \Lag'$ intertwining the action of $C^+(m)$.
Moreover, $\psi$ may be taken orthogonal with respect to the inner
product on $\vo$. Indeed as $\psi$ is non singular and $\psi^*$ is
also in $End_{C^+(m)}(\vo_m^p)$ we have that $\xi = (\psi
\psi^*)^{-1/2} \psi \in End_{C^+(m)}(\vo_m^p)$ is orthogonal and
$\xi \Lag =\Lag'$. Since $\vo = \Lag \oplus J_m(\Lag) = \Lag'
\oplus J_m(\Lag')$ we can extend $\psi$ to all of $\vo$ by
$\psi(J_m (u)) = J_m(\psi (u))$ for all $u \in \Lag$. Since $\psi$
is orthogonal by construction and $\psi J_i = J_i \psi $ for all
$i=1,\ldots,m$, $\psi$ is automorphism.

If, instead, $m \equiv 0  \,(\mod 4)$, $\Lag$ and $\Lag'$ are
isomorphic as $C^+(m)$-modules if and only if the multiplicity of
$\vo_m^+$ is the same in both. Notice also that if $\psi$ is an
orthogonal automorphism  of $\n$, $\psi$ :$\Lag \to \Lag'$ is a
 $C^+(m)$-module isomorphism. The rest of the proof follows as in the
 previous case, to conclude that there are exactly $p+1$ isomorphism
 clases.
\end{proof}

The following theorem essentially solves  our problem for $\n$
irreducible as algebra of Heisenberg type.

\begin{theorem} \label{irreducible} If  $\vo$ is an irreducible
$C(m)$-module, then every proper $C^+(m)$-submodule of $\vo$ is
Lagrangian.
\end{theorem}
\begin{proof}
Let $z_1,\ldots,z_m$ be an orthonormal basis of $\z_m$ and denote
$K = J_{z_1}\ldots J_{z_m}$

For $m \equiv 3,5,6,7 \,(\mod 8)$ we see, by dimension, that there
are no proper $C^+(m)$-modules (see Table \ref{tabla1}).

For $m \equiv 0,4  \,(\mod 8)$, $K$ is symmetric and $K^2=1$, and
by Table \ref{tabla1} the only proper $C^+(m)$-modules are
irreducible. In fact, they are the $K$-eigenspaces $\vo_\pm$ with
eigenvalues $\pm 1$. Since $J_{z_i}K = -KJ_{z_i}$, $J_{z_i}\vo_\pm
= \vo_\mp$. Let $u,v \in \vo_+$. Then $[u,v]=0$ if and only if
$\langle z_{i}, [u,v] \rangle =0$
 ($i=1,\ldots,m$). Now, $\langle z_i, [u,v] \rangle = \langle
J_{z_i}u,v \rangle =0$ because  $\langle \vo_+, \vo_- \rangle=0$.
So, $\vo_+$ is Lagrangian. In analogous way we obtain that $\vo_-$
is Lagrangian as well.

The cases $m \equiv 1,2 \,(\mod 8)$ are a bit more involved.
 We
will find and fix a Lagrangian $\Lag$ that, by Corollary
\ref{cmasmodulo} will necessarily be a $C^+(m)$-module. We will
compute the space $End_{C^+(m)}(\vo)$ of intertwining operators
for the $C^+(m)$-action on $\vo$ and prove that $\phi(\Lag)$ is
isotropic for all $\phi \in End_{C^+(m)}(\vo)$.

For $m =1$ any one dimensional subspace can be taken as our fixed
Lagrangian. For $m=2$,  we can take $\vo= \Hc$, $\z_2 = {\rm
span}_\R\{i,j\}$ and the Clifford action given by left
multiplication and inner product  $\langle u,v\rangle=
u\overline{v}$. It is easy to see that $\Lag = {\rm
span}_\R\{i,j\}$ is a Lagrangian in $\vo$. Furthermore, if
$\phi:\Hc \to \Hc$ is an intertwining operator and $\phi(1)=q$,
then $\phi=R_q$ is right multiplication by $q$. It follows that
$\phi(\Lag)$ is Lagrangian. Let us choose $\phi=R_q$ with
$q^2=-1$, so that $\phi$ is skew-symmetric.

When $m \equiv 1,2 \,(\mod 8)$ there is only one irreducible
$C(m)$-module $\vo_m$ up to equivalence. If $m=r+8$, then $\vo_m =
\vo_8\otimes \vo_r$, by Corollary \ref{tensorirreducible}. Let
$\Oc_\pm$ be the vector subspace of $\vo_8$ of eigenvectors of
$K_8$ with eigenvalues $\pm1$. $\Oc_\pm$ are $C^+(8)$-modules
since $K_8$ commute with $C^+(8)$. If $\Lag_r\in {\rm
Lag}(\vo_r)$, then $\Lag = \Oc_+ \otimes \Lag_r \oplus \Oc_-
\otimes \Lag_r^\perp\in {\rm Lag}(\vo_m)$ and $\vo_m = \Lag \oplus
\Lag^\perp$. The existence of a Lagrangian for any dimension
follows now by induction. We fix such $\Lag$ and we will see that
$\phi(\Lag)$ is isotropic for all $\phi \in End_{C^+}(\vo_m)$.

In the case $m \equiv 1 \,(\mod 8)$, we first prove that
$$
End_{C^+}(\vo_m)=\left\{ \begin{pmatrix}a\,\, Id_{\Lag} &
cK_{|\Lag^\perp}
\\ bK_{|\Lag} & d\,\, Id_{\Lag^\perp}
\end{pmatrix}: a,b,c,d \in \R\right\},
$$
where the action of the matrix is with respect to the
decomposition $\vo_m = \Lag \oplus \Lag^\perp$.  Now, $K: \Lag \to
\Lag^\perp$ is an intertwining operator for the action of
$C^+(m)$. Thus, the result follows by Schur's Lemma and the fact
that $\Lag$ and $\Lag^\perp$ are real irreducible $C^+(m)$-modules
(see  Table \ref{tabla1}). Therefore any proper $C^+(m)$-module is
of the form
$$ \Lag'= \{ax+bKx: x \in
\Lag\}.$$ Such $\Lag'$ is isotropic: $[ax+bKx,ay+bKy] =
ab([x,Ky]+[Kx,y] )$, but for all $z \in \R^m$
$$
\langle z, [x,Ky]+[Kx,y]\rangle = \langle J_zx,Ky\rangle + \langle
J_zKx,y\rangle =  \langle -KJ_zx,y\rangle + \langle J_zKx,y\rangle
=0
$$
since $KJ_z = J_zK$ for all $z \in \R^m$.

In the case $m \equiv 2 \,(\mod 8)$, we have $K^2=-1$. Since $m$
is even and $\Lag$ and $\Lag^\perp$  are $C^+(m)$-modules,  $\Lag$
must be $K$-invariant.  Thus, $K$ gives
 complex structures on $\Lag$ and $\Lag^\perp$. By Schur's
Lemma, the intertwining operators of $\Lag$  and $\Lag^\perp$  are
of the form $a + b K$, with $a,b \in \R$. Now, we look for
$\phi_m$ an $C(m)$-intertwining operator of $\vo_m$ that sends
$\Lag$ to $\Lag^\perp$. Since $\Lag = \Oc_+ \otimes \Lag_r \oplus
\Oc_- \otimes \Lag_r^\perp$, $\phi_m$ can be defined recursively
as
\begin{equation}\label{phisubeme}
\phi_2 = R_q, \quad \text{and} \quad \phi_m=Id \otimes \phi_r
\end{equation}
for $m=8+r$ and $r \equiv 2 \,(\mod 8)$. Since  $\phi_2$ is
skew-symmetric, $\phi_m$ is so too.

With $\phi=\phi_m$ we have
$$
End_{C^+}(\vo)=\left\{ \begin{pmatrix}
 a + b K & (c +d K)\phi_{|\Lag^\perp}
 \\ (a'+b'K)\phi_{|\Lag} & c' + d'K
\end{pmatrix}: a,b,c,d,a',b',c',d' \in \R\right\}
$$
with respect to the decomposition $\vo_m = \Lag \oplus
\Lag^\perp$. Recall that $KJ_i=-J_iK$ ($i=1,\ldots,m$) and
$K^t=-K$. Let $\mathcal W = \{ u_x: u_x=( a + b K)x + (a'+b'K)\phi
x, x \in \Lag\}$, then
\begin{align*}
[u_x,u_y] &= [( a + b K)x,(a'+b'K)\phi y]+ [(a'+b'K)\phi x,( a + b
K)y].
\end{align*}
For $z \in \R^m$, we have
\begin{align*}
\langle z, [u_x,u_y]\rangle  &= \langle z,[( a + b
K)x,(a'+b'K)\phi y]+ [(a'+b'K)\phi x,( a + b K)y]\rangle  \\
& =\langle J_z( a + bK)x,(a'+b'K)\phi y\rangle +\langle
J_z(a'+b'K)\phi=
 x,( a + b K)y\rangle  \\
& =\langle (a'-b'K)J_z( a + bK)x,\phi y\rangle -\langle ( a - b=
 K)J_z(a'+b'K) x,\phi y\rangle  \\
& =\langle J_z(a'+b'K)( a + bK)x,\phi y\rangle -\langle J_z( a +
b=
 K)(a'+b'K) x,\phi y\rangle  \\
& = 0
\end{align*}
where we use that $\phi_m$ is skew-symmetric. From the equation
above we conclude that $\mathcal W$ is isotropic.
\end{proof}

\begin{corollary} If $m \equiv 1 \,(\mod 8)$ and $m=r+8$, then every
Lagrangian of $\vo_{m}$ is of the form
$$(\Oc_+ \otimes \Lag') \oplus (\Oc_- \otimes
\Lag'^\perp)$$ where $\Lag'$ is a Lagrangian of $\vo_r$.
\end{corollary}
\begin{proof}
>From the proof above we know that $\Lag =\{ax+bK_m x : x \in \Oc_+
\otimes \Lag_r \oplus \Oc_- \otimes \Lag_r^\perp \}$. Now $x = x_+
\otimes y + x_- \otimes y^\perp$ and $K_m = Id \otimes K_r$. Then
$$
ax+bK_m x = x_+ \otimes (a+bK_r)y + x_- \otimes (a+bK_r)y^\perp.
$$
By the Theorem above, $\Lag'=\{(a+bK_r)y: y \in \Lag_r\}$ is a
Lagrangian in $\vo_r$ and is easy to check that
$\Lag'^\perp=\{(a+bK_r)y^\perp: y^\perp \in \Lag_r^\perp\}$.
\end{proof}

\end{section}

\begin{section}{${\rm Lag}(\n)$ explicitly}

In this section we will describe ${\rm Lag}(\n)$. We do this in
terms of Plücker coordinates on the corresponding grassmanian
whenever the result is not too messy, and as orbits of
automorphisms of $\n$, or finite unions thereof, in every case.
The notation ${\rm Lag}(\n)\cong A/B$ will always mean that ${\rm
Lag}(\n)$ is an orbit of a subgroup $A\in {\rm Aut}_o(\n)$ and $B$
is the isotropy subgroup.

Thanks to Corollary \ref{cmasmodulo}, it is enough to describe
 the $C^+(m)$-submodules of $\vo$ and then determine which of
 these are totally isotropic.

We will consider eight separate cases, corresponding to the
congruence modulo 8 of  $\dim \z$, the center of $\n$. Fix
$z_1,\ldots,z_m$  an orthonormal basis of $\z$, set $J_i =
J_{z_i}$ and
$$K_m=J_1 ... J_m.$$
It will be convenient to display the dimension $m$ of the center
of $\n=\vo+\z$ and the multiplicities of the spin representations
appearing in $\vo$, as explained in section 2:
$$\n_m^{(p)} =(\vo_m)^p \oplus \z_m,$$
for $m\not \equiv 3 \,(\mod 4)$, and
$$\n_m^{(p_+,p_-)} = (\vo_m^+)^{p_+} \oplus (\vo_m^-)^{p_-} \oplus \z_m$$
for $m \equiv 3 \,(\mod 4)$.

\

\centerline{\bf The case $m \equiv 1$}

\

In this case there is only one irreducible $C(m)$-module $\vo_m$,
up to isomorphism.  From the proof of Theorem \ref{irreducible} we
have $\vo_m =\Lag \oplus K\Lag$ where $\Lag$ and $K\Lag$ are
Lagrangian and equivalent irreducible $C^+(m)$-modules. If $\vo$
is a $C(m)$-module, then
$$
\vo= \bigoplus _{i=1}^p (\wo_i \oplus  K\wo_i),
$$
where $\wo_i$ is isomorphic to $\Lag$ for all $i$.  Thus, we can
write
$$\vo= (\R^p \otimes \Lag) \oplus (\R^p \otimes K\Lag)= \R^p
\otimes \vo_m.$$ In terms of this decomposition,
$$[a\otimes
w,a'\otimes w'] = \langle a,a' \rangle [w,w'],$$  for $a,a' \in
\R^p$, $w,w' \in \vo_m$. From Schur's Lemma and the fact that
$\Lag$ and $K\Lag$ are real irreducible $C^+(m)$-modules, it
follows that
$$
End_{C^+}(\vo)=\left\{ \begin{pmatrix}A \otimes Id_\Lag & C
\otimes K_{|K\Lag}
\\ B \otimes K_{|\Lag} & D \otimes Id_{K\Lag}
\end{pmatrix}: A,B,C,D \in M_p(\R)\right\},
$$
where the action  is with respect to the decomposition $(\R^p
\otimes \Lag) \oplus (\R^p \otimes K\Lag)$. Therefore any proper
$C^+(m)$-module is of the form
$$\mathcal W= {\rm
span}_\R\{Aa\otimes u +Ba \otimes Ku: a\otimes u \in \R^p \otimes
\Lag\}.$$

\begin{lemma}
\begin{enumerate}
\item \label{lem:1}  $\mathcal W$ is isotropic if and only if $A^tB-B^tA=0
$.
\item \label{lem:2}$\mathcal W$ is Lagrangian if and only if it is isotropic
 and
the $2p \times p$-block matrix $\begin{pmatrix}A \\
B\end{pmatrix}$ has rank equal to $p$
\end{enumerate}
\end{lemma}
\begin{proof} (1) Let $a\otimes u, a'\otimes u'\in  \R^p \otimes
\Lag$. Then
\begin{align*}
[Aa\otimes u +Ba \otimes Ku,Aa'\otimes u' +Ba' \otimes Ku'] &=
\langle Aa,Ba'\rangle [u,Ku']+\langle Ba,Aa'\rangle[Ku,u'] \\
&= (\langle Aa,Ba'\rangle -\langle Ba,Aa'\rangle)[u,Ku']\\
&=\langle (A^tB-B^tA)a,a'\rangle [u,Ku'].
\end{align*}
Since $\Lag, K\Lag \in {\rm Lag}(\n)$, $[u,Ku']\not=0$, for some
$u,u' \in \Lag$. Therefore, $\mathcal W$ is isotropic if and only
if $A^tB-B^tA=0$ and this proves (\ref{lem:1}).

(\ref{lem:2}) This follows because $\dim(\mathcal W)$ is equal to
\,$\text{rank}\hspace{-.1cm}\begin{pmatrix}A \\
B\end{pmatrix} \dim(\Lag)$.
\end{proof}

\begin{proposition}\label{prop:5.2}
${\rm Aut}_o(\n_m^{(p)})$ is the group of
$$\begin{pmatrix}A \otimes Id_\Lag & C
\otimes K_{|K\Lag}
\\ B \otimes K_{|\Lag} & D \otimes Id_{K\Lag}
\end{pmatrix} \in End_{C^+(m)}(\vo)$$ such that
\begin{equation}
A^tB-B^tA=0,  \qquad C^tD-D^tC=0, \qquad A^tD-B^tC=1.
\end{equation}
Thus, ${\rm Aut}_o(\n_m^{(p)})\cong Sp(p,\R)$ (see also \cite{S}).
\end{proposition}
\begin{proof}

 By (\ref{eq:2.3}), $g  \in {\rm Aut}_o(\n)$ if and only if
$g$ is in $End_{C^+}(\vo)$ and $g^tJ_mg=J_m$. We can choose a
basis of $\vo$ such that $J_m= \begin{pmatrix}0 & Id \otimes
K_{|K\Lag}
\\ - Id \otimes K_{|\Lag} & 0 \end{pmatrix}$, thus the result follows.
\end{proof}

\begin{proposition}
$${\rm Lag}(\n_m^{(p)})
\cong U(p)/O(p).$$
\end{proposition}
\begin{proof}
As  a maximal compact subgroup of  ${\rm Aut}_o(\n)= Sp(p,\R)$,
the group of isometric automorphisms $\mathcal A(\n)$ is
isomorphic to $U(p)$, viewed  as the set of matrices
$\begin{pmatrix}A \otimes Id_\Lag & -B \otimes K_{|K\Lag}
\\ B \otimes K_{|\Lag} & A \otimes Id_{K\Lag}
\end{pmatrix}\in End_{C^+}(\vo)$. By Proposition \ref{transitively},
$U(p)$ acts transitively on ${\rm Lag}(\n)$. The isotropy subgroup
is the set of matrices   $\begin{pmatrix}A \otimes Id_\Lag & 0
\\  0 & A \otimes Id_{K\Lag}
\end{pmatrix}$ where $A$ is real orthogonal matrix, hence this
isotropy subgroup is isomorphic to $O(p)$.
\end{proof}

\

\centerline{\bf The case $m \equiv 2$}

\

In this case there is only one irreducible $C(m)$-module $\vo_m$,
up to isomorphism. By Theorem \ref{irreducible}, $\vo_m = \Lag
\oplus \phi_m \Lag$ where $\Lag$ is an irreducible $C^+(m)$-module
and $\phi_m$ is given by (\ref{phisubeme}). Then any $C(m)$-module
has a decomposition
$$
\vo= \bigoplus _{i=1}^p (\wo_i \oplus  \phi_m\wo_i),
$$
where $\wo_i$ is isomorphic to $\Lag$ for all $i$.  In this case
$K$ is a complex structure on $\vo$ which leaves $\Lag$ invariant,
we can write
$$\vo= (\C^p
\otimes \Lag) \oplus (\C^p \otimes \phi_m\Lag) = \C^p \otimes_\C
\vo_m.$$
In terms of this decomposition,
 $$[a\otimes w,a'\otimes w'] =
[\langle a,a' \rangle_\C\, w,w'],$$ where $a,a' \in \C^p$, $w,w'
\in \vo_m$ and $\langle a,a' \rangle_\C = \sum a_i{a'_i}$. Indeed,
We can write $a\otimes w= (b +Kc) \otimes w$ and $a'\otimes w'=
(b' +Kc') \otimes w'$, where $a,a' \in \C^p$, $b,b',c,c' \in
\R^p$. Thus,
\begin{align*}
[a\otimes w,a'\otimes w'] &= [(b +Kc) \otimes w,(b' +Kc') \otimes
w'] \\
&= [b\otimes w,b'\otimes w']+ [Kc\otimes w,Kc'\otimes w']+
[b\otimes w,Kc'\otimes w']+ [Kc\otimes w,b'\otimes w'] \\
&= [b\otimes w,b'\otimes w']+ [c\otimes Kw,c'\otimes Kw']+
[b\otimes w,c'\otimes Kw']+ [c\otimes Kw,b'\otimes w'] \\
&= [b\otimes w,b'\otimes w']+ [c\otimes K^2 w,c'\otimes w']+
[b\otimes K w,c'\otimes w']+ [c\otimes Kw,b'\otimes w'] \\
&=[(\langle b, b' \rangle - \langle c, c' \rangle)w,w']+
[(\langle c, b' \rangle + \langle b, c' \rangle) Kw,w'] \\
&= [\langle a,a' \rangle_\C\, w,w'].
\end{align*}

It follows from \cite{S}, Schur's Lemma and the fact that $\Lag$
and $\phi_m\Lag$ are complex irreducible $C^+(m)$-modules, that
$$
End_{C^+}(\vo)=\left\{ \begin{pmatrix}A \otimes Id_\Lag & C
\otimes {\phi_m}_{|\phi_m\Lag}
\\ B \otimes {\phi_m}_{|\Lag} & D \otimes Id_{\phi_m\Lag}
\end{pmatrix}: A,B,C,D \in M_p(\C)\right\},
$$
where the action of the matrix is with respect to the
decomposition $(\C^p \otimes \Lag) \oplus (\C^p \otimes
\phi_m\Lag)$ (cf. Table \ref{tabla1}). Therefore any proper
$C^+(m)$-module is of the form
$$ \mathcal W= {\rm span}_\C\{Aa\otimes u
+Ba \otimes \phi_mu: a\otimes u \in \C^p \otimes \Lag\}.$$ Now we
can check which $\mathcal W$ are isotropic:
\begin{align*}
[Aa\otimes u +Ba \otimes \phi_mu,Aa'\otimes u' +Ba' \otimes
\phi_mu'] &= [ \langle Aa,Ba'\rangle u,\phi_mu']+[\langle
Ba,Aa'\rangle \phi_mu,u']
\\
&= [(\langle Aa,Ba'\rangle -\langle Ba,Aa'\rangle)u,\phi_mu'] \\
&=[\langle (B^tA -A^tB)a,a'\rangle u,\phi_mu'].
\end{align*}
 Therefore, $\mathcal W$ is isotropic if and only if $A^tB-B^tA=0$.

\begin{proposition}\label{prop:5.5}
 ${\rm Aut}_o(\n_m^{(p)})$ is the
group of matrices
$$\begin{pmatrix}A \otimes Id_\Lag & C \otimes
{\phi_m}_{|\phi_m\Lag}
\\ -B \otimes {\phi_m}_{|\Lag} & D \otimes Id_{\phi_m\Lag}
\end{pmatrix} \in End_{C^+}(\vo)$$ such that
\begin{equation}
A^tB-B^tA=0, \qquad C^tD-D^tC=0, \qquad A^tD-B^tC=1.
\end{equation}
Thus, ${\rm Aut}_o(\n_m^{(p)})\cong Sp(p,\C)$.
\end{proposition}

\begin{proof}
By (\ref{eq:2.3}), $\xi  \in {\rm Aut}_o(\n)$ if and only if $\xi$
is in $End_{C^+}(\vo)$ and $\xi^tJ_m\xi=J_m$. As in the proof of
Theorem \ref{irreducible}, we set $m = 8+r$, $\vo_m = \vo_8
\otimes \vo_r$, $\phi_m = Id \otimes \phi_r$ and $\Lag_m =\Oc_+
\otimes \Lag_r \oplus \Oc_- \otimes \Lag_r^\perp$, where $\Lag_r$
is a Lagrangian of $\vo_r$. We also have $\phi_m(\Lag_m) =
\Lag_m^\perp= J_m(\Lag_m)$. Now, we fix basis $v_1,\ldots,v_8$ of
$\Oc_+$ and $w_1,\ldots,w_t$ of $\Lag_r$, thus $v_i \otimes w_j$
is a basis of $\Oc_+ \otimes \Lag_r$. We complete to a basis of
$\Lag_m$ adding $(J_8 \otimes J_r)(v_i \otimes w_j)$ which is a
basis of $\Oc_- \otimes \Lag_r^\perp$. Finally, applying $J_m= K_8
\otimes J_r$ we complete to a basis of $\vo_m$. With respect of
this basis the matrix of $J_m =
\begin{pmatrix}  0 & -Id \\ Id & 0\end{pmatrix}$ and $\phi_m =
\begin{pmatrix}  0 & M_m \\ N_m & 0\end{pmatrix}$. If $m=2$,
then $\vo_2= \Hc$, $J_2$ is left multiplication by $j$ and $\Lag_2
= {\rm span}_\R\{i,j\}$. We can take $\phi_2=R_j$ the right
multiplication by $j$. In this case $M_2=
\begin{pmatrix}  1 & 0 \\ 0 & -1\end{pmatrix}$ and $N_2 = -M_2$.
As   $\phi_m = Id \otimes \phi_r$, an easy recursion shows that
$M_m$ is symmetric, $M_m^2=Id$  and $N_m = -M_m$, for all $m
\equiv 2 \,(\mod 8)$. Now, $\vo_m^p = \C^p \otimes \vo_m$ and $\xi
\in End_{C^+}(\vo_m^p)$ can be written as
$$
\xi = \begin{pmatrix}A \otimes Id_\Lag & C \otimes M_m
\\  B \otimes M_m & D \otimes Id_{\Lag^\perp}
\end{pmatrix},  \text{ so $\xi \in {\rm Aut}_o(\n)$ if and only if }
 \begin{pmatrix}A  & C
\\ B & D \end{pmatrix} \text{ is in } Sp(p,\C).
$$

\end{proof}

\begin{proposition}
$${\rm Lag}(\n_m^{(p)})\cong U(p,\Hc)/U(p).$$
\end{proposition}
 \begin{proof}
 By Proposition \ref{transitively} the group of orthogonal automorphism,
 which is $U(p,\Hc)$, acts transitively. Recall that $U(p,\Hc)$
 can be identified with the subgroup of $ Sp(p,\C)$ given by the matrices
 $\begin{pmatrix} A & -\overline{B} \\ B & \overline{A} \end{pmatrix}$.  To
 find the isotropy group of $\Lag_m$, we have that   $\begin{pmatrix} A &=
 -\overline{B} \\ B & \overline{A} \end{pmatrix} \begin{pmatrix} \Lag_m \\=
 0\end{pmatrix} = \begin{pmatrix} \Lag_m \\ 0\end{pmatrix}$ implies $B=0=
$ and $A^t \overline{A}=1$, thus $A \in U(p) $.
 \end{proof}

\

\centerline{\bf The case $m \equiv 3$}

 \

We begin with a

\begin{lemma}\label{lemma:61} Let $m= 8s +3$ and let $\vo_8^{\otimes s}=
 \otimes \vo_3^\pm$ be the
$C(m)$-module obtained by using repeatedly  Proposition
\ref{prop:tensorial}. Then, $\vo_m^{\pm} = \vo_8^{\otimes s}
\otimes \vo_3^\pm$. Moreover, let $j_1,j_2,j_3$ be canonical
generators of $C(3)$, be $j'_1,\ldots,j'_8$ be canonical
generators of $C(8)$ and define
\begin{align*}
J_{8k+t} &= K_8^{\otimes k} \otimes j'_t \otimes 1^{\otimes (s-k)}
, \quad=
 \text{ if } 0 \le k<s, 1\le t \le 8; \\
J_{8s+t} &=  K_8^{\otimes s} \otimes j_t, \quad \text{ if } 1\le t
\le 3;
\end{align*}
where $K_8 = j'_1\ldots j'_8$. Then, $J_1, \ldots, J_m$ are
generators of $C(m)$ such that
$$
J_i^2 = -1 , \quad J_iJ_k = -J_kJ_i, \quad \text{ for } 1 \le i,k
\le m \text{ and } i\not=k.
$$
Moreover,  $K_m = J_1 \ldots J_m = Id^{\otimes s} \otimes K_3$.
\end{lemma}
\begin{proof}$\vo_8^{\otimes s} \otimes \vo_3^\pm$ is irreducible and $K_m =
= Id^{\otimes s} \otimes K_3$ by Corollary
\ref{tensorirreducible}. As ${K_m}_{|\vo_m^{\pm}}= \pm Id$, we
have  $\vo_m^{\pm} = \vo_8^{\otimes s} \otimes \vo_3^\pm$. The
relations on the $J_i$'s follow by straightforward computation.

\end{proof}

\begin{proposition} \label{pigualq} If ${\rm Lag}(\n_m^{(p_+,p_-)})
\not=\emptyset$, then $p_+=p_-$.
\end{proposition}
\begin{proof}
The trace of the operator $K_m$ on $\vo$ is $tr(K_m) =
(p_+-p_-)\dim \vo_m^+$. We will prove that there exists a
Lagrangian $\Lag$ if and only if  $tr(K_m)=0$. Suppose that $\Lag$
is a Lagrangian, then $\Lag^\perp = J_i(\Lag)$ is also a
Lagrangian and $\vo = \Lag \oplus \Lag^\perp$. Since $K_m$ is an
odd product of $J_i$'s, $K_m$ sends $\Lag$ to $\Lag^\perp$ and
$\Lag^\perp$ to $\Lag$, thus $tr(K_m)=0$, so $p_+=p_-$.
\end{proof}

We will see later that the converse also holds. From now on, we
consider only the case $p_+=p_-$.

Next we will describe  the intertwining operators between
$(\vo_m^\pm)^p$ and $(\vo_m^\pm)^p$ as $C^+(m)$-modules. Using the
explicit construction of an algebra of Heisenberg type with center
of dimension $3$ given by (\ref{prodint3}), we let $\phi: \vo_3^+
\to \vo_3^-$ be given by $\phi(u) = \overline{u}$. Then $\phi$
intertwines $\vo_3^+$ and $\vo_3^-$ as $C(3)^+$-modules.

\begin{proposition} \label{prop:5.8}
\begin{enumerate}
\item\label{iiuno} The $C^+(m)$-intertwining operators of $(\vo_m^+)^p$
(resp. $(\vo_m^-)^p$) are $Id \otimes R_{\overline{A}}$ (resp. $Id
\otimes L_A$) with $A \in gl(p,\Hc)$ and $R_{\overline{A}}:
(\vo_3^+)^p \to (\vo_3^+)^p$ (resp. $L_A: (\vo_3^-)^p \to
(\vo_3^-)^p$) denotes the right (resp. left) action, i.e. $R_{X} u
= u X^t$ (resp. $L_Xu = Xu$) with $X \in gl(p,\Hc)$.
\item \label{iidos} Let $\varphi = K_8^{\otimes s} \otimes \phi$, then
 $\varphi: \vo_m^+ \to \vo_m^-$
intertwines the action of $C^+(m)$ and  the intertwining operators
between $(\vo_m^+)^p$ and $(\vo_m^-)^p$ (resp.  $(\vo_m^-)^p$  and
$(\vo_m^+)^p$) are $(Id \otimes L_A)\varphi $ (resp. $(Id \otimes
R_{\overline{A}})\varphi$) with $A \in gl(p,\Hc)$.
\end{enumerate}

\end{proposition}
\begin{proof}
(\ref{iiuno}) Let $m = 8s +3$. For $s=0$, the result follows by
Schur's Lemma. For $s=1$, let $\theta$ be an intertwining operator
of  $(\vo_m^+)^p$. We can write $\theta \sum_i \alpha_i \otimes
\beta_i$ where $\beta_i:(\vo_3^+)^p \to (\vo_3^+)^p$ are linear
independent and $\alpha_i: \vo_8^p \to \vo_8^p$. Now, $\theta$
commutes with $C^+(8) \otimes Id$, therefore the $\alpha_i$'s
commute with the action of $C^+(8)$. Thus, $\alpha_i = a_i Id +
b_i K_8$ and
$$\theta = Id \otimes (\sum_i a_i\beta_i) + K_8
\otimes (\sum b_i \beta_i).$$ Since $\theta$ commutes with $Id
\otimes C^+(3)$, $\sum_i a_i\beta_i=R_q$  and $\sum b_i
\beta_i=R_{q'}$  for some $q,q' \in \mathbb H$. Finally, since
$K_8$ anti-commutes with $j_i'$ (cf. Lemma \ref{lemma:61}) $q'$
must be zero.  The case $s>1$ is similar.

(\ref{iidos}) It easy to see that $\varphi$ commutes with the
action of $C^+(m)$  and, therefore, that any intertwining operator
is composition of $\varphi$ with an operator from the first part
of the proposition.

\end{proof}

Set
 $$\vo_m^{(p_+,p_-)} := (\vo_m^+)^{p_+} \oplus
(\vo_m^-)^{p_-}.$$
 Recall
that if $m=8s+3$,
\begin{equation}\label{eq:descomp}
\vo_m^{(p_+,p_-)} = (\vo_8^{\otimes s} \otimes \vo_3^+)^{p_+}
\oplus (\vo_8^{\otimes s} \otimes \vo_3^-)^{p_-}.
\end{equation}

\begin{corollary}\label{intertwiningm3}
\begin{equation}\label{for53}
End_{C^+(m)}(\vo_m ^{(p,p)}) = \{\begin{pmatrix} Id \otimes
 R_{\overline{A}} & (Id \otimes R_{\overline{C}})\varphi\\
(Id \otimes L_B)\varphi &Id \otimes L_D \end{pmatrix}: \text{ with
} A,B,C, D \in gl(p,\Hc)\},
\end{equation}
where the blocks are with respect to the decomposition
(\ref{eq:descomp}).
\end{corollary}

\begin{remark}\label{transpuestas} If we view $(\vo_3^+)^p$ as a real space
 with inner product defined
as in (\ref{prodint3}), then for $X \in gl(p,\Hc)$ the transpose
of $R_X:(\vo_3^+)^p \to (\vo_3^+)^p$ is $R_{X^*}$, where $X^* =
\overline{X}^t$. In analogous way, the transpose of $L_X:
(\vo_3^-)^p \to (\vo_3^-)^p$ is $L_{X^*}$. Writing $\phi$ in a
canonical basis, we see that it is symmetric.
\end{remark}

\begin{lemma}\label{transpra} For all $A,B,C,D \in gl(p,\Hc)$,
\begin{enumerate}
\item $\phi L_A = R_{\overline A_{}}\phi$ and $\phi R_A = L_{\overline
A}\phi$,
\item
$
 \begin{pmatrix}Id \otimes  R_{\overline{A}} & (Id \otimes R_{\overline{C}})=
 \varphi
 \\ (Id \otimes L_B) \varphi &Id \otimes  L_D \end{pmatrix}^{{}_t} =
 \begin{pmatrix} Id \otimes R_{A^t} & (Id \otimes R_{B^t}) \varphi \\ (Id
 \otimes L_{C^*}) \varphi &Id \otimes  L_{D^*}
 \end{pmatrix}.
$
\end{enumerate}
\end{lemma}
\begin{proof} (1) This follows from  $\overline{Ax} = \overline{x}=
 \overline{A}^t$ and
$\overline{x A} =  \overline{A}^t \overline{x}$.

(2) If $m=3$, we have
$$
\begin{pmatrix} R_{\overline{A}} & R_{\overline{C}} \phi \\ L_B \phi &L_D=
 \end{pmatrix}^{t} =
 \begin{pmatrix} (R_{\overline{A}})^{t} &  (L_B \phi)^{t} \\=
 (R_{\overline{C}} \phi)^{t}  &(L_D)^{t}
 \end{pmatrix}.
$$
Now, using Remark \ref{transpuestas} and (1), we have $(L_B
\phi)^{t} = \phi^t(L_B)^{t} = \phi L_{B^*} = R_{B^t} \phi$. In
analogous way we have $(R_{\overline{C}} \phi)^{t} = L_{C^*}\phi$.
For $m > 3$ the result follows easily from the case $m=3$ and the
fact that $\varphi=K_8^{\otimes s} \otimes \phi$ is symmetric.
\end{proof}

\begin{remark}
$End_{C^+(3)}(\vo_3^{(p,p)})$ is an associative algebra isomorphic
to $gl(\Hc,2p)$, under the  isomorphism $\Theta$:
\begin{gather*}
\begin{pmatrix} R_{\overline{A}} & R_{\overline{C}} \phi \\
L_B \phi &L_D \end{pmatrix} \mapsto  \begin{pmatrix} \phi & 0 \\
0 & Id \end{pmatrix}\begin{pmatrix} R_{\overline{A}} &
R_{\overline{C}} \phi
 \\
L_B \phi &L_D \end{pmatrix}  \begin{pmatrix} \phi & 0 \\
0 & Id \end{pmatrix} = \begin{pmatrix} L_{{A}} & L_{{C}} \\L_B
&L_D \end{pmatrix} \mapsto   \begin{pmatrix} {A} & {C}
\\B &D
\end{pmatrix}.
\end{gather*}
Moreover, if $U=\begin{pmatrix} R_{\overline{A}} &
R_{\overline{C}} \phi
 \\
L_B \phi &L_D \end{pmatrix}$, then $\Theta(U^t) =\Theta(U)^* $.
From (\ref{for53}), it is clear that $End_{C^+(m)}(\vo_m^{(2p)})$
is an associative algebra isomorphic to
$End_{C^+(3)}(\vo_3^{(p,p)})$, so we can view $\Theta$ as an
isomorphism  $\Theta: End_{C^+(m)}(\vo_m^{(p,p)}) \to gl(\Hc,2p)$.
\end{remark}

Recall that the group $Sp(p,q)$ is defined as
$$
Sp(p,q) = \{X \in GL(p+q,\Hc): X^* I_{p,q} X = I_{p,q}\},
$$
where $I_{p,q}= \begin{pmatrix} Id_p & 0   \\ 0& -Id_q
\end{pmatrix}$ (see \cite{K}, p. 70).

\begin{proposition}[\cite{S}] \label{autcaso7} $\Theta$ defines an
 isomorphism
$${\rm Aut}_o(\n_m^{(p,p)})\cong Sp(p,p).$$ More explicitly,
$$\begin{pmatrix} Id \otimes R_{\overline{A}} & (Id \otimes
 R_{\overline{C}}) \varphi
\\ (Id \otimes L_B) \varphi & Id \otimes L_D
\end{pmatrix}\ \in\ {\rm Aut}_o(\n_m^{(2p)})$$ if and only if
 \begin{align}
A^*A -  B^*B&= Id, \\
 {C}^*C - D^*D  &= -Id \\
A^*C - B^*D &=0.
\end{align}
\end{proposition}
\begin{proof}
By (\ref{eq:2.3})
$$
{\rm Aut}_o(\n_m^{(p,p)}) = \{U \in End_{C^+(m)}(\vo_m^{(2p)}):
U^{t} J_i U = J_i, \text{ for some } i=1,\ldots,m \}.
$$
Moreover, since $m$ is odd, this is equivalent to
\begin{equation}\label{autcerito}
{\rm Aut}_o(\n_m^{(p,p)}) = \{U \in End_{C^+(m)}(\vo_m^{(p,p)}):
U^{t} K_m U = K_m \}.
\end{equation}
Consider the real canonical basis  $1,i,j,k$  of $\vo_3^\pm =\Hc$
and $\{e_i\}$ be a basis of $\vo_8^{\otimes s}$, then
$\{e_i\otimes 1,e_i\otimes i,e_i\otimes j,e_i\otimes k\}$ is a
basis of $\vo_m^\pm$. With respect to this basis $ K_m:
\vo_m^{(p,p)} \to
 \vo_m^{(p,p)}$ has a matrix $\begin{pmatrix} Id_{\vo_m^+} & 0 \\
0  & -Id_{\vo_m^-} \end{pmatrix}$. Now, $U \in {\rm
Aut}_o(\n_m^{(2p)})$ if and only if $U^t K_m U = K_m$ if and only
if $\Theta(U^t K_m U) = \Theta(K_m)$ if and only if $\Theta(U)^*
I_{p,p} \Theta(U) = I_{p,p}$ if and only if $\Theta(U) \in
Sp(p,p)$.
\end{proof}

Let us recall by Proposition \ref{prop:tensorial} that the bracket
in $\vo_m^\pm= \vo_8^{\otimes s} \otimes \vo_3^\pm$ is given by:
$$
[v_1\otimes u_1,v_2 \otimes u_2] = [v_1,v_2]\langle u_1,u_2
\rangle + \langle K_8^{\otimes s} v_1, v_2 \rangle [u_1,u_2],
$$
where $v_1,v_2 \in \vo_8^{\otimes s}$ and $u_1,u_2 \in \vo_3^\pm$.

\begin{lemma} \label{lagmitad}
$$\mathcal L = {\rm span}_{\R}\{ v \otimes u +\varphi(v \otimes
u):\ \ v \in \vo_8^{\otimes s},  u \in (\vo_3^+)^p\}$$ is a
Lagrangian subspace.
\end{lemma}
\begin{proof}
 Recall that $\varphi = K_8^{\otimes s}  \otimes
\phi$. Let $v_1,v_2 \in \vo_8^{\otimes s}$ and $u_1,u_2 \in
(\vo_3^+)^p$, then we have to see that

\begin{equation}\label{eq:lagm3}
[v_1\otimes u_1+K_8^{\otimes s} v_1 \otimes \phi u_1,v_2 \otimes
u_2+K_8^{\otimes s} v_2 \otimes \phi u_2] = 0.
\end{equation}

First, using the fact that $K_8^2=1$, $K_8^t =K_8$ it is easy to
see, by induction, that  $[K_8^{\otimes s}v_1,K_8^{\otimes s}v_2]
= -[v_1,v_2]$ and $\langle v_1,K_8^{\otimes s}v_2\rangle = \langle
K_8^{\otimes s}v_1,v_2\rangle$. Also $[\phi u_1, \phi u_2]_{\vo^-}
= -[u_1,u_2]_{\vo^+}$ and $\langle \phi u_1, \phi u_2 \rangle =
\langle u_1,
 u_2 \rangle$ (see (\ref{prodint3}) and (\ref{corchete3})). Thus,
\begin{align*}
[v_1\otimes &u_1+K_8^{\otimes s} v_1 \otimes \phi u_1,v_2 \otimes
u_2+K_8^{\otimes s} v_2 \otimes \phi u_2] = \\
&=[v_1\otimes u_1 ,v_2 \otimes u_2] +[K_8^{\otimes s} v_1 \otimes
\phi u_1,K_8^{\otimes s} v_2 \otimes \phi u_2] \\
&=[v_1\otimes u_1 ,v_2 \otimes u_2]+ [K_8^{\otimes
s}v_1,K_8^{\otimes s}v_2]\langle \phi u_1,\phi u_2 \rangle +
\langle
(K_8^2)^{\otimes s} v_1, K_8^{\otimes s}v_2 \rangle [\phi u_1,\phi u_2] \\
&=[v_1\otimes u_1 ,v_2 \otimes u_2]- [v_1,v_2]\langle u_1,u_2
\rangle - \langle  v_1, K_8^{\otimes s}v_2 \rangle [u_1,u_2] \\
 &=[v_1\otimes u_1 ,v_2 \otimes u_2]- [v_1,v_2]\langle u_1,u_2
\rangle - \langle  K_8^{\otimes s} v_1, v_2 \rangle [u_1,u_2] \\
&= 0,
\end{align*}
\end{proof}

\begin{proposition} \label{prop:5.15}Let $A,B,C,D \in gl(p,\Hc)$, then
$$\mathcal W =\{\begin{pmatrix} Id \otimes R_{\overline{A}} & (Id \otimes
 R_{\overline{C}}) \varphi \\
(Id \otimes L_B) \varphi & Id \otimes L_D
\end{pmatrix}\begin{pmatrix} v \otimes u \\ \varphi(v \otimes u)=
 \end{pmatrix}: v \in
\vo_8^{\otimes s},  u \in (\vo_3^+)^p\},
$$
is Lagrangian if and only if
\begin{align} \label{eq:m3lag}
(A+C)^*(A+C)-({B+D})^*(B+D)=0.
\end{align}
 and $A+C$ and $B+D$ are non singular.
\end{proposition}
\begin{proof}
Every $C^+(m)$-module is of the form:
$$
\mathcal W =\{\begin{pmatrix} Id \otimes R_{\overline{A}} & (Id
\otimes
 R_{\overline{C}}) \varphi \\
(Id \otimes L_B) \varphi & Id \otimes L_D
\end{pmatrix}\begin{pmatrix} v \otimes u \\ \varphi(v \otimes u)=
 \end{pmatrix}: v \in
\vo_8^{\otimes s},  u \in (\vo_3^+)^p\},
$$
for some $A,B,C,D \in gl(p,\Hc)$. Now we impose that $\mathcal W$
 be isotropic. First,
\begin{equation*}
\begin{pmatrix} Id \otimes R_{\overline{A}} & (Id \otimes R_{\overline{C}})=
 \varphi \\
(Id \otimes L_B) \varphi & Id \otimes L_D
\end{pmatrix}\begin{pmatrix} v \otimes u \\ \varphi(v \otimes u)
\end{pmatrix}=
\begin{pmatrix} v \otimes R_{\overline{A+C}}u \\ (Id \otimes=
 L_{B+D})\varphi(v \otimes u)
\end{pmatrix},
\end{equation*}
Then the bracket of two elements of $\mathcal W$ is of the form:
\begin{align*}\label{eq:lag} \tag {*}
[v_1 \otimes &R_{\overline{A+C}}u_1 +(Id \otimes
L_{B+D})\varphi(v_1 \otimes u_1), v_2 \otimes
R_{\overline{A+C}}u_2 + (Id \otimes L_{B+D})\varphi(v_2 \otimes
u_2)] \\
&= [v_1 \otimes R_{\overline{A+C}}u_1 , v_2 \otimes
R_{\overline{A+C}}u_2 ] + [(Id \otimes L_{B+D})\varphi(v_1 \otimes
u_1), (Id \otimes L_{B+D})\varphi(v_2 \otimes u_2)].
\end{align*}
Now,
\begin{align*}
[v_1 \otimes R_{\overline{A+C}}u_1 , v_2 \otimes
R_{\overline{A+C}}u_2 ] &=[v_1,v_2]\langle R_{\overline{A+C}}u_1,
R_{\overline{A+C}}u_2\rangle + \langle K_8^{\otimes s}v_1,v_2
\rangle [ R_{\overline{A+C}}u_1, R_{\overline{A+C}}u_2],
\end{align*}
and
\begin{align*}
[&(Id \otimes L_{B+D})\varphi(v_1 \otimes u_1), (Id \otimes
L_{B+D})\varphi(v_2 \otimes u_2)]  \\ &= [K_8^{\otimes
s}v_1,K_8^{\otimes s}v_2]\langle L_{B+D}\phi u_1,  L_{B+D}\phi
u_2\rangle + \langle (K_8^2)^{\otimes s}v_1, K_8^{\otimes s}v_2
\rangle [L_{B+D}\phi u_1,  L_{B+D}\phi u_2] \\
 &= -[v_1,v_2]\langle L_{B+D}\phi u_1,  L_{B+D}\phi
u_2\rangle + \langle K_8^{\otimes s}v_1, v_2
\rangle [L_{B+D}\phi u_1,  L_{B+D}\phi u_2] \\
&= -[v_1,v_2]\langle \phi R_{\overline{B+D}}u_1,  \phi
R_{\overline{B+D}} u_2\rangle + \langle K_8^{\otimes s}v_1, v_2
\rangle [\phi R_{\overline{B+D}}u_1,  \phi R_{\overline{B+D}} u_2]
\\
&= -[v_1,v_2]\langle  R_{\overline{B+D}}u_1,  R_{\overline{B+D}}
u_2\rangle - \langle K_8^{\otimes s}v_1, v_2 \rangle [
R_{\overline{B+D}}u_1,  R_{\overline{B+D}} u_2].
\\
\end{align*}
The RHS of (\ref{eq:lag}) is equal to
\begin{align*}
[v_1,v_2](\langle &R_{\overline{A+C}}u_1,
R_{\overline{A+C}}u_2\rangle -\langle  R_{\overline{B+D}}u_1,
R_{\overline{B+D}} u_2\rangle) \\ &+  \langle K_8^{\otimes
s}v_1,v_2 \rangle ([ R_{\overline{A+C}}u_1, R_{\overline{A+C}}u_2]
- [ R_{\overline{B+D}}u_1, R_{\overline{B+D}} u_2]).
\end{align*}
Varying $v_1,v_2$ we have that the the RHS of (\ref{eq:lag}) is
equal to 0 if and only if:
\begin{align}
\langle R_{\overline{A+C}}u_1, R_{\overline{A+C}}u_2\rangle
-\langle  R_{\overline{B+D}}u_1, R_{\overline{B+D}} u_2\rangle &=
0\\
[ R_{\overline{A+C}}u_1, R_{\overline{A+C}}u_2] - [
R_{\overline{B+D}}u_1, R_{\overline{B+D}} u_2] &= 0,
\end{align}
for all $u_1,u_2 \in (\vo_3^+)^p$. Thus, equations
(\ref{prodint3}) and (\ref{corchete3}) imply the result.
\end{proof}

\begin{proposition} \label{prop:5.16}
$${\rm Lag}(\n_m^{(p,p)})\cong
U(p,\Hc) \times U(p,\Hc) /(U(p,\Hc) \times Id) = U(p,\Hc)$$
$${\rm Lag}(\n_m^{(p,q)})= \emptyset\qquad p\not= q$$
\end{proposition}
\begin{proof}
If $\Lag$ is Lagrangian of $\vo$, then every element $x$ in $\Lag$
is of the form
\begin{align*}
x = \begin{pmatrix} Id \otimes R_{\overline{A}} & (Id \otimes
 R_{\overline{C}}) \varphi \\
(Id \otimes L_B) \varphi & Id \otimes L_D
\end{pmatrix}\begin{pmatrix} v \otimes u \\ \varphi(v \otimes u)
\end{pmatrix}
&=\begin{pmatrix} v \otimes R_{\overline{A+C}}u \\ (Id \otimes
L_{B+D})\varphi(v \otimes u) \end{pmatrix} \\
&= \begin{pmatrix} v \otimes u' \\ (Id \otimes
L_{(B+D)(A+C)^{-1}}\varphi(v \otimes u')
\end{pmatrix} \\
&=\psi \begin{pmatrix} v \otimes u' \\ \varphi(v \otimes u')
\end{pmatrix},
\end{align*}
where $\psi =
\begin{pmatrix} Id \otimes Id & 0 \\
0 & Id \otimes L_{D'} \end{pmatrix} $, with $D' =
(B+D)(A+C)^{-1}$. Since $D'$ is unitary by equation
(\ref{eq:m3lag}), $\psi  \in \mathcal A(\n_m^{(p,p)})\cong
U(p,\Hc) \times U(p,\Hc)$. It is clear that the isotropy group is
$U(p,\Hc) \times Id$.
\end{proof}

\

\centerline{\bf The case $m \equiv 4$} \label{sec:mequiv4}

\

Let $\vo_m = \wo_+ \oplus \wo_-$ be the decomposition of the
$C(m)$-module into the eigenspaces of $K_m$ of eigenvalues $\pm
1$. Thus, $\vo_m^{p} = \wo_+^{p} \oplus \wo_-^{p}$.

\begin{lemma} \label{above} Let $\Lag_1 \subset \wo_+^p$ be a
 $C^+(m)$-submodule
and $\Lag_1^\perp$  its orthogonal complement in $\wo_+^p$. Then,
\begin{enumerate}
\item  $\Lag = \Lag_1 + J_m (\Lag_1^\perp)$ is a Lagrangian
subspace of $\vo_m^{p}$.
\item Every Lagrangian is of this form.
\end{enumerate}
\end{lemma}
\begin{proof}  (1) It is clear that the dimension of $\Lag$  is $\frac12
\dim(\vo_m^{p})$.  As $\wo_+^p$ and $\wo_-^{p}$ are isotropic, we
must only verify that $[\Lag_1 , J_m (\Lag_1^\perp)]=0$. This
follows since for every $z \in \z_m$ we have:
\begin{align*}
\langle z ,[\Lag_1 , J_m (\Lag_1^\perp)]\rangle = \langle
J_z(\Lag_1) , J_m (\Lag_1^\perp)\rangle = - \langle J_mJ_z(\Lag_1)
, \Lag_1^\perp\rangle = - \langle \Lag_1 , \Lag_1^\perp\rangle =0.
\end{align*}

\noindent (2) Let $\Lag$ be any Lagrangian,  by Corollary
\ref{cmasmodulo}, $\Lag$ is a $C^+(m)$-module and it can be
decomposed as $\Lag= \Lag_1 \oplus \Lag_{-1}$ as eigenspaces of
$K_m$. Now,
$$
\langle J_m(\Lag_1) ,\Lag_{-1}\rangle = \langle z_m, [\Lag_1
,\Lag_{-1}]\rangle =0.
$$
Thus, $\Lag_{-1} \subset J_m(\Lag_1)^\perp = J_m(\Lag_1^\perp) $,
therefore , by dimension, $\Lag_{-1} = J_m(\Lag_1^\perp)$.
\end{proof}

\begin{proposition} \label{prop:5.18} ${\rm Aut}_o(\n_m^{(p)})$ has $p+1$ orbits
 in
${\rm Lag}(\n_m^{(p)})$, of the form
$$U(p,\Hc) /
U(r,\Hc)\times U(p-r,\Hc)$$ $r=0,\ldots,p$.
\end{proposition}
\begin{proof} From the Lemma \ref{above}, every Lagrangian is
determined by a $C^+(m)$-module of $\wo_+^{p}$. Any $\psi \in {\rm
Aut}_o(\n_m^{(p)})$ preserves $\wo_{\pm}$. Furthermore, given any
pair $\Lag_1$, $\Lag_1'$ of $C^+(m)$-submodules of $\wo_{+}^{p}$
of the same dimension there exists a non singular $\psi \in
End_{C^+(m)}(\wo_+^{p})$ such that $\psi \Lag_1 =\Lag_1'$.
Moreover, $\psi$ may be taken orthogonal with respect to the inner
product. Indeed as $\psi$ is non singular and $\psi^*$ is also in
$End_{C^+(m)}(\wo_+^p)$ we have that $\xi = (\psi \psi^*)^{-1/2}
\psi \in End_{C^+(m)}(\wo_+^p)$ is orthogonal and $\xi \Lag_1
=\Lag_1'$. Now we extend $\xi$ to an element of $\mathcal
A(\n_m^{(p)})$ as $\xi(J_m(w)) = J_m\xi(w)$ for all $w \in
\wo_+^{p}$. For each $i=0,\dots,p$, fix $\Lag_1^{(r)}$ any
$C^+(m)$-submodule of $\wo_+^{p}$ of dimension $r . \dim(\wo_+)$.
Then  $\Lag^{(r)} = \Lag_1^{(r)} + J_m ((\Lag_1^{(r)})^\perp)$,
$r=0,\dots,p$, are representatives of each orbit. From \cite{S},
any $\xi \in {\rm Aut}_o(\n_m^{(p)})$ can be written as
$\begin{pmatrix} A & 0
\\ 0 & (A^*)^{-1} \end{pmatrix}$, $A
\in Gl(p,\Hc)$, with respect to any  basis compatible with the
decomposition   $\wo_+^{p} \oplus \wo_-^{p}$. Thus, with this
identification,   $\mathcal A(\n_m^{(p)}) \cong U(p,\Hc) \subset
Gl(p,\Hc)$ and $U(r,\Hc)\times U(p-r,\Hc)$ is the isotropy
subgroup of $\Lag^{(r)}$.
\end{proof}

\

\centerline{\bf The case $m \equiv 5$}

\

Consider the inclusion $C(m-1) \hookrightarrow C(m)$ via $J_i
\mapsto J_i$ ($i=1,\ldots,m-1$). Then $\vo_m$ is an irreducible
$C(m-1)$-module, so $\vo_m^p = \vo_{m-1}^p$. Denote by
$\wo_{\pm}^p \subset \vo_m^p$ the eigenspace of $K_{m-1}$ of
eigenvalue $\pm1$. Clearly, $\wo_{\pm}^p$ is a $C^+(m-1)$-module.

\begin{lemma} \label{lemma:5.19} $\Lag$  is a Lagrangian of $\n_m^{(p)}$ if
 and only if
there exists a $C^+(m-1)$-submodule  $\Lag_+\subset \wo_{+}^p$
such that $J_m(\Lag_+) = \Lag_+^\perp $,  $\dim(\Lag_+) = \frac12
\dim \wo_{+}^p$ and $\Lag = \Lag_+\oplus J_{m-1}J_m(\Lag_+)$. Here
$\Lag_+^\perp$ is the orthogonal complement of $\Lag_+$ in
$\wo_{+}^p$.
\end{lemma}
\begin{proof}
Any Lagrangian $\Lag$ in $\vo_m^p$ is also a Lagrangian in
$\vo_{m-1}^p$ and we have seen in section \ref{sec:mequiv4}, that
$\Lag = \Lag_+ \oplus J_{m-1} (\Lag_+^\perp)$, where $\Lag_+$ is a
$C^+(m-1)$-submodule of $\wo_{+}^p$. Now, $J_m$ commutes with
$K_{m-1}$, so $J_m$ preserves $\wo_{\pm}^p$ and $J_m(\Lag_+)
\subset \wo_{+}^p$. Since $\Lag_+$ is isotropic, $\langle
J_m(\Lag_+), \Lag_+ \rangle =0$, thus $J_m(\Lag_+)= \Lag_+^\perp$
and
 $\wo_{+}^p = \Lag_+ \oplus J_m(\Lag_+)$. Conversely, given $\Lag=
 \Lag_+\oplus J_{m-1}(\Lag_+^\perp)$,
 it is clear that $\langle z_m, [\Lag,\Lag] \rangle = 0$,  and this
 implies that $\Lag$ is Lagrangian.
\end{proof}

\begin{remark}
As $\dim(\Lag_+) = \frac12 \dim \wo_{+}^p$ and $\Lag_+$ is a
$C^+(m)$-module, $p$ is even and $\Lag_+$ is isomorphic to
$\wo_+^{p/2}$ as $C^+(m)$-module.
\end{remark}

We see that there are no lagrangians in $\n_m^{(p)}$ unless $p$ is
even. Put
$$p=2q$$
so $q$ is arbitrary natural number.

\begin{lemma}\label{lemma:5.22} $\vo_{m-1}^q \otimes \vo_1$ is isomorphic to
$\vo_{m}^{2p}$ as $C(m)$-modules.
\end{lemma}
\begin{proof} $\vo_{m-1}^q \otimes \vo_1$ has an structure of $C(m)$-module
given by Proposition \ref{prop:tensorial}. For $q=1$ the results
follow by dimension. For $q>1$, $\vo_{m-1}^q \otimes \vo_1=
(\vo_{m-1} \otimes \vo_1)^q = \vo_{m}^{p}$.
\end{proof}

By  Lemma \ref{lemma:5.22}, we can consider  $\n_m = \vo_{m-1}^q
\otimes \vo_1 \oplus \z$. Therefore, by Proposition
\ref{prop:tensorial}, $J_m: \vo_{m-1}^q \otimes \vo_1 \to
\vo_{m-1}^q \otimes \vo_1$ is given by $K_{m-1} \otimes j_1$.

\begin{proposition} Let $1,i$ a basis of $\vo_1$
\begin{enumerate}
\item $\wo_{+}^{p} = \wo_{+}^q \otimes \vo_1$ is the eigenspace of
 $K_{m-1}$
of eigenvalue $1$.
\item  With $\Lag_+ =  \wo_{+}^q \otimes \mathbb R 1$,
$$\Lag =
\Lag_+\oplus J_{m-1}J_m(\Lag_+)$$ is a Lagrangian of $\n$.
\end{enumerate}
\end{proposition}
\begin{proof}
{ (1)} The result follows since $K_{m-1}: \vo_{m-1}^q \otimes
\vo_1 \to \vo_{m-1}^q \otimes \vo_1$ acts trivially on $\vo_1$.

{(2)} Clearly $\Lag_+$ is a $C^+(m-1)$-module and $J_m = K_{m-1}
\otimes j_1$ sends   $\Lag_+$ to  $(\Lag_+)^\perp =  \wo_{+}^q
\otimes \mathbb R i$. Thus by Lemma \ref{lemma:5.19} $\Lag =
\Lag_+\oplus J_{m-1}J_m(\Lag_+)$ is a Lagrangian.
\end{proof}

\begin{proposition}
$$ {\rm Lag}(\n_m^{(2q)}) \cong U(2q)/U(q,\Hc)$$
$$ {\rm Lag}(\n_m^{(2q+1)}) = \emptyset$$
\end{proposition}

\begin{proof}By Schur we have that $End_{C^+(m-1)}(\wo_+^{2q})$ is equal to
$Gl(p,\Hc)$. As $\wo_+^{2q} = (\wo_{+}^q \otimes \mathbb R 1)
\oplus (\wo_{+}^q \otimes \mathbb R i)$ as $C^+(m-1)$-module,
every  $C^+(m-1)$-intertwining operator is of the form
$\xi=\begin{pmatrix} A & C \\ B & D
\end{pmatrix}$ where $A,B,C,D \in Gl(q,\Hc)$. Moreover, with
respect to this decomposition $J_m = \begin{pmatrix} 0 & -Id \\ Id
& 0 \end{pmatrix}$. If $\Lag_1$ be any Lagrangian of $\n$, then
$\Lag_1 = \Lag_1^+ \oplus J_{m-1}J_m (\Lag_1^+)$, where $\Lag_1^+
= \{  \begin{pmatrix} A & C \\ B & D \end{pmatrix}
\begin{pmatrix} w \\ 0 \end{pmatrix}: w \in \Lag_+\}$.
We can choose $B = -C$, $D = A$ and $A^*A +C^*C =1$. Now,
$J_m(\Lag_1^+) = (\Lag_1^+)^\perp$ if and only if for $w,v \in
\Lag_+$ we have
\begin{align*}
0 &= \langle J_m \begin{pmatrix} A & C \\ -C & A=
 \end{pmatrix}\begin{pmatrix} w \\ 0 \end{pmatrix}, \begin{pmatrix} A & C \\=
 -C & A \end{pmatrix}\begin{pmatrix} v \\ 0 \end{pmatrix} \rangle \\
&=  \langle J_m \begin{pmatrix} Aw \\ -Cw \end{pmatrix},
\begin{pmatrix}=
 Av \\ -Cv \end{pmatrix} \rangle \\
&= \langle  \begin{pmatrix} Cw \\ Aw \end{pmatrix},
\begin{pmatrix} Av \\ -Cv \end{pmatrix} \rangle  = \langle (A^*C
-C^*A)w, v\rangle.
\end{align*}
Thus,  $J_m(\Lag_1^+) = (\Lag_1^+)^\perp$ if and only if $A^*C
-C^*A =0$.

On the other hand, we can extend $\xi \in End_{C^+(m-1)}(\wo_+^p)$
to an operator $\tilde{\xi}$ in $ End_{C^+(m)}(\vo_m^{2q})$ acting
on $\wo_-^{2q}$ by $v \mapsto J_{m-1}J_m \xi J_mJ_{m-1}v$. It is
clear that $\tilde{\xi}$ commutes with $C^+(m)$. By
(\ref{eq:2.3}), $\tilde{\xi}$ is in ${\rm Aut}_o(\n)$ if and only
if
$$J_m = \xi^t J_m \xi =   \begin{pmatrix}
A^* & -C^* \\ C^* & A^* \end{pmatrix}\begin{pmatrix} 0 & -1 \\ 1 &
0
\end{pmatrix}\begin{pmatrix} A & C \\ -C & A \end{pmatrix}.
$$
So $\tilde{\xi}$ is in ${\rm Aut}_o(\n)$ if and only if $\xi$
satisfies $ A^*C =C^*A$ and $A^*A + C^*C =1$. The isotropy
subgroup of $\Lag$ is given by $C=0$ and thus $A \in U(q,\Hc)$.
Since, ${\rm Aut}_o(\n_m^n) \cong Gl({n},\Hc) \cap O(2n,\C)$ and
the group of orthogonal automorphisms is $U(n)$ (cf.
\cite{R}\cite{S}), we have that the variety of Lagrangians is
indeed $U(2q)/U(q,\Hc)$, and any Lagrangian is in ${\rm
Aut}_o(\n)\Lag$.
\end{proof}

\

\centerline{\bf The case $m \equiv 6$}

\

The inclusions $C(m-2) \subset C(m-1) \subset  C(m)$ show that
$\vo_m$ is an irreducible $C(m-2)$-module -- hence $\vo_m^p =
\vo_{m-1}^p=\vo_{m-2}^p$. Denote by $\wo_{\pm}^p \subset \vo_m^p$
the eigenspace of $K_{m-2}$ of eigenvalue $\pm1$. As in the
previous case we have that any Lagrangian can be written as
$$\Lag
= \Lag_+ \oplus J_{m-1}J_{m-2}(\Lag_+)$$ where $\Lag_+$ is a
$C^+(m-2)$-module.

\begin{proposition}
$\Lag$ is Lagrangian of $\vo_m^p$  if and only if $\Lag$ is a
$K_m$-invariant Lagrangian of $\vo_{m-1}^p$.
\end{proposition}
\begin{proof} Any Lagrangian of
$\vo_m^p$ is Lagrangian of $\vo_{m-1}^p$ and as it is
$C^+(m)$-module is $K_m$-invariant. Conversely, let $\Lag$ be a
$K_m$-invariant Lagrangian of $\vo_{m-1}^p$. By Lemma
\ref{lemma:5.19}, $\Lag = \Lag_+ \oplus J_{m-2}J_{m-1}(\Lag_+)$,
where $\Lag_+$ is a $C^+(m-2)$-module and $J_{m-1}(\Lag_+) =
(\Lag_+)^\perp$. Let's see first  that $J_{m}(\Lag_+) =
(\Lag_+)^\perp$. Indeed, $K_m(\Lag) = \Lag$ and $K_m(\wo_{\pm}^p)
= \wo_{\pm}^p$, thus $K_m(\Lag_+) = \Lag_+$. On the other hand,
$K_{m-2}(\Lag_+) = \Lag_+$, $J_{m-1}(\Lag_+) = (\Lag_+)^\perp$ and
$J_mJ_{m-1} =  K_m K_{m-2}$ leaves $\Lag_+$ invariant. So
$J_{m}(\Lag_+) = (\Lag_+)^\perp$. It is follows immediately that
$[\Lag_+,\Lag_+]=0$ and $[ J_{m-2}J_{m-1}(\Lag_+),
J_{m-2}J_{m-1}(\Lag_+)]=0$. As $J_{m}(\Lag_+) = (\Lag_+)^\perp
\subset \wo_+^p$ and $\wo_+^p \perp \wo_-^p$ we have $\langle
J_m(\Lag_+), J_{m-2}J_{m-1}(\Lag_+)\rangle=0$. This implies that
 $[\Lag_+, J_{m-2}J_{m-1}(\Lag_+)]=0$.
\end{proof}

  \begin{proposition}
  $${\rm Lag}(\n_m^{(2q)})\cong O(2q)/O(q)
\times O(q)$$

$${\rm Lag}(\n_m^{(2q+1)}) = \emptyset$$
\end{proposition}

\begin{proof} To  compute ${\rm Aut}_o(\n)$, note that any
irreducible $C^+(m)$-module  has a complex structure given by
$K_m$. So, $End_{C^+(m)}(\vo)$ is isomorphic to $\C_p$. By
(\ref{eq:2.3}), ${\rm Aut}_o(\n) \cong O(p,\C)=\{\xi \in \C_p:
\xi^t\xi = Id\}$. Also, $\mathcal A(\n) \cong O(p)$. With respect
to the decomposition
$$\vo_m^p = \Lag \oplus J_{m}(\Lag)$$ we write $\xi
\in O(p)$ as $\begin{pmatrix} A & B \\ C & D \end{pmatrix}$. Then
the isotropy group of $\Lag$ is constituted by the matrices
$\begin{pmatrix} A & 0
\\ 0 & D \end{pmatrix}$, with $A,D \in O(q)$, thus the result
follows.
\end{proof}

\

\centerline{\bf The case $m \equiv 7$}

\

The first part is identical to the case $\dim(\z) \equiv 3$. Write
  $\vo_m^{\pm} = \vo_8^{\otimes s} \otimes \vo_7^\pm$, let
$\phi: \vo_7^+ \to \vo_7^-$ be given by $\phi(x) = \overline{x}$,
where conjugation is
 octonionic, and define  $\varphi= K_8^{\otimes s}
\otimes \phi:\vo_m^{+} \to \vo_m^{-}$. We also see that
$$\n_m^{(p_+,p_-)} =
((\vo_m^+)^{p_+ } \oplus (\vo_m^-)^{p_-} \oplus \z_m$$ has  a
Lagrangian if and only if $p_+= p_-$. Write  as before
$$\vo_m^{(p,p)} = (\vo_m^+)^p \oplus (\vo_m^-)^p, \qquad
\n_m^{(p,p)} = \vo_m ^{(2p)} \oplus \z_m.$$ Then
\begin{proposition}
\begin{enumerate}
\item
$$
End_{C^+(m)}(\vo_m^{(p,p)}) \simeq \{\begin{pmatrix} Id \otimes A
& (Id
 \otimes C) \varphi\\
(Id \otimes B)\varphi &Id \otimes D \end{pmatrix}: \text{ with }
A,B,C, D \in gl(p,\mathbb R)\},
$$
where the matrices are written with respect to the above
decomposition.
\item  ${\rm Aut}_o(\n_m^{(p,p)}) \cong O(p,p)$ and $\mathcal A(\n_m ^p)
\cong O(p)
 \times O(p)$.
\item $\mathcal L = {\rm span}_{\R}\{ v \otimes u +\varphi(v \otimes u): v \in
\vo_8^{\otimes s},  u \in (\vo_7^+)^p\}$ is a Lagrangian subspace.
\item The group $\mathcal A(\n_m^{(p,p)})$ acts transitively on the variety
 of
Lagrangians subspaces.
\end{enumerate}
\end{proposition}
\begin{proof}
(1)  is proved as in  Proposition \ref{pigualq}. For (2) we notice
first that the space of intertwining operators of $(\vo_7^\pm)^p$
is  $gl(p,\mathbb R)$. The rest of the proof follows the lines of
Proposition \ref{prop:5.8}. (3) follows as in Proposition
\ref{autcaso7}, noticing that $X^*$ must be replaced by $X^t$. (4)
is proved exactly as Lemma \ref{lagmitad}. (5) is proved along the
lines of the propositions \ref{prop:5.15} and \ref{prop:5.16}.
\end{proof}

\begin{corollary}
$${\rm Lag}(\n_m^{(p,p)}) \cong O(p) \times O(p) / (O(p) \times Id)
\cong O(p)$$

$${\rm Lag}(\n_m^{(p,q)}) = \emptyset \qquad (p\not=q)$$

\end{corollary}

\

\centerline{\bf The case $m \equiv 8$} \label{sec:mequiv8}

\

The argument parallels that of the case $m \equiv 4$, except that
the last Proposition should be replaced by

\begin{proposition}   ${\rm Aut}_o(\n_m^{(p)})$ has $p+1$ orbits in
${\rm Lag}(\n_m^{(p)})$, of the form
$$O(p) / O(r)\times
O(p-r)$$ $r=0,\ldots,p$
\end{proposition}
\begin{proof} From the Lemma \ref{above}, every Lagrangian is
determined by a $C^+(m)$-module of $\wo_+^{p}$. Any $\psi \in {\rm
Aut}_o(\n_m^{(p)})$ preserves $\wo_{\pm}$. Furthermore, given any
pair $\Lag_1$, $\Lag_1'$ of $C^+(m)$-submodules of $\wo_{+}^{p}$
of the same dimension there exists a non singular $\psi \in
End_{C^+(m)}(\wo_+^{p})$ such that $\psi \Lag_1 =\Lag_1'$.
Moreover, $\psi$ may be taken orthogonal with respect to the inner
product. Indeed as $\psi$ is non singular and $\psi^*$ is also in
$End_{C^+(m)}(\wo_+^p)$ we have that $\xi = (\psi \psi^*)^{-1/2}
\psi \in End_{C^+(m)}(\wo_+^p)$ is orthogonal and $\xi \Lag_1
=\Lag_1'$. Now we extend $\xi$ to an element of $\mathcal
A(\n_m^{(p)})$ as $\xi(J_m(w)) = J_m\xi(w)$ for all $w \in
\wo_+^{p}$. For each $i=0,\dots,p$, fix $\Lag_1^{(r)}$ any
$C^+(m)$-submodule of $\wo_+^{p}$ of dimension $r . \dim(\wo_+)$.
Then  $\Lag^{(r)} = \Lag_1^{(r)} + J_m ((\Lag_1^{(r)})^\perp)$,
$r=0,\dots,p$, are representatives of each orbit. From \cite{S},
 any
$\xi \in {\rm Aut}_o(\n_m^{(p)})$ can be written as
$\begin{pmatrix} A & 0
\\ 0 & (A^*)^{-1} \end{pmatrix}$ where $A \in  Gl(p,\R)$,
 with respect to any  basis compatible with the
decomposition   $\wo_+^{p} \oplus \wo_-^{p}$. With this
identification $\mathcal A(\n_m^{(p)}) \cong O(p) \subset
Gl(p,\R)$ and $O(r)\times O(p-r)$ is the isotropy subgroup of
$\Lag^{(r)}$.

\end{proof}

\end{section}

\end{document}